\documentclass[a4paper,12pt,leqno]{article}
\usepackage{exscale}
\usepackage{amsmath}
\usepackage{amsfonts}
\usepackage{wasysym}
\usepackage{stmaryrd}
\usepackage{amscd}
\usepackage{amssymb}
\usepackage{theorem}
\usepackage{latexsym}
\usepackage[final]{epsfig}
\usepackage{epsf}
\usepackage{graphicx}

\newfont{\abc}{cmtt10 scaled 1200}

\def\R{\mathbb{R}}
\def\Z{\mathbb{Z}}

\def\ve{\varepsilon}

\def\ra{\rightarrow}

\def\p{\partial}
\def\qed{\hfill $\Box$ \\}

\def\op{\operatorname}
\def\Ric{\op{Ric}}

\def\ch{\mathcal{H}}

\def\op{\operatorname}
\def\Ric{\op{Ric}}

\def\supp{\op{supp}}

\def\op{\operatorname}
\def\Ric{\op{Ric}}
\def\scal{\op{scal}}

\def\supp{\op{supp}}
\def\x1n-1{x_1,\dots,x_{n-1}}

\def\dipsi1{\frac{\partial\psi}{\partial x_1}}

\def\d2phin{\frac{\partial^2\psi}{\partial x_n^2}}

\def\1d{\frac{1}{d}}

\def\inn{\overset{\circ}}

\def\op{\operatorname}
\def\Ric{\op{Ric}}
\def\scal{\op{scal}}

\def\supp{\op{supp}}
\def\x1n-1{x_1,\dots,x_{n-1}}

\def\dipsi1{\frac{\partial\psi}{\partial x_1}}

\def\d2phin{\frac{\partial^2\psi}{\partial x_n^2}}

\def\1d{\frac{1}{d}}

\def\inn{\overset{\circ}}
\def\lev{{}^\approx |A|^{-1}}
\def\ap{{}^\approx |A|}

\usepackage{amscd}

\newtheorem{proposition}{Proposition}[subsection]
{\theorembodyfont{\rmfamily}\newtheorem{definition}[proposition]{Definition}}
\newtheorem{lemma}[proposition]{Lemma}
{\theorembodyfont{\rmfamily}\newtheorem{remark}[proposition]{Remark}}

\newtheorem{corollary}[proposition]{Corollary}

{\theorembodyfont{\rmfamily}}

\begin{document}

\vspace*{0cm}
\begin{center}\Large{\bf{Inductive Analysis on Singular Minimal Hypersurfaces}}\\
\medskip
{\small{by}}\\
\medskip
\large{\bf{ Joachim Lohkamp}} \\

\end{center}

\vspace{1cm}
\noindent Mathematisches Institut, Universit\"at M\"unster, Einsteinstrasse 62, Germany\\
 {\small{\emph{e-mail: lohkampj@math.uni-muenster.de}}}
\vspace{1cm}

\setcounter{section}{1}
\renewcommand{\thesubsection}{\thesection}
\subsection{Introduction}

\normalsize
\bigskip

The geometric analysis of a minimal hypersurface $H$ within some Riemannian manifold $(M,g)$ with second fundamental form $A$ usually involves the scalar quantity  $| A |^2$ = sum of
squared principal curvatures. A few classical examples are seen from Simons type inequalities like: $\Delta_H|A|^2 \ge -C \cdot (1+|A|^2)^2$ or the stability condition (valid in particular
for area minimizers):\\ $0 \le Area'' (f) = \int_H | \nabla_H f |^2 - f^2 (| A |^2 + Ric_M (\nu, \nu)) d A$ for infinitesimal variations $f$ in normal direction $\nu$ to $H$ .\\ Now the point is
that minimality, i.e. $trA = 0$, also implies that the scalar curvature $scal_H$ of $H$ satisfies $scal_H = - | A |^2$ in a flat ambient space (which is  central in particular when one
is interested in the case where $H$ is singular).  \\

Thus one realizes that the analysis of the scalar curvature of $H$ is intimately linked to the analysis of the underlying minimal hypersurface $H$. Oftentimes, a look at the \emph{conformal
Laplacian} $Lu = - \triangle u + \frac{n-2}{4 (n-1)} \cdot scal_H \cdot u$ is a good starting point to get a feeling of global aspects of $scal_H$ since it allows us, in particular, to measure
averages of $scal_H$. \\ $L$ is well-understood when $H$ is compact and smooth, cf. [KW]. Namely, recall that in this case the first eigenfunction $f_H$ of $L$ does not vanish and
(choosing it to stay positive) it can be used to conformally deform the metric on $H$ $g_H$ into a metric $f_H^{4/n-2} \cdot g_H$ of scalar curvature of fixed sign (equal to that of the
eigenvalue) as is readily seen from the $scal$-transformation law under conformal changes:
\[\lambda_1 \cdot f_H = - \triangle  f_H + \frac{n-2}{4 (n-1)} \cdot scal_{(H,g_H)} \cdot f_H = scal_{(H,f_H^{4/n-2} \cdot g_H)} \cdot f_H^{n+2/n-2}\]

However for  scalar curvature geometry smooth minimal hypersurfaces are a too narrow class of objects: \emph{singular} minimal hypersurfaces appear as intermediate objects even if the focus
is on smooth manifolds. This brings us to the main topic of this paper settling the basic classical question how to extract  information
 from singular hypersurfaces encoded in its conformal Laplacian in a fashion that fits with the smooth case. This can directly be used resp. translated to understand the way how such minimal hypersurfaces
 inherit positive scalar curvature from their ambience, cf. [CL], resp. how to smooth singular minimal hypersurfaces to regular hypersurfaces with positive mean curvature, cf. [L1].\\

To survey the paper we recall that a major aid to deal with the uncontrollable singular set $\Sigma \subset H$ is the inductive use of \emph{tangent cones} around points in $\Sigma$ (cf. [Gi],
[Si]). These cones generalize tangent planes of smooth submanifolds: after scaling around some $p \in \Sigma$ (by some $\tau_m \ra \infty$) one approaches a (usually non-unique) limit object
which is a (locally) area minimizing cone $C_p$ which approximates $\tau_m \cdot H$. The intersection of $C_p$ with the distance sphere $\p B_1(0)$ around the tip of the cone produces
again a minimal hypersurface $F$ within the sphere $\p B_1(0)$. Of course, since $C_p$ can also have a singular set $\sigma \varsupsetneqq \{0\}$, $F$ may have singularities and
(although $F$ is not area minimizing) the region close to the singular set can be handle in the same way as for area minimizers (namely the tangent cones are area minimizing cones). This
permits us to argue by (dimensional) induction tackling the singular set of
$F$ the same way as $\Sigma \subset H$ . Eventually one reaches (in dimension 8 or possibly earlier) a \emph{regular} cone (= singular only in the tip) where one has explicit control.\\

The natural question is whether one could find a better control over the singular set $\Sigma$. However in a way this seems to be the wrong question. It is a classical result (cf. [D],[Gi] and [Si])
\emph{compact} set has at least \emph{codimension 7} within $H^n$.  In dimensions $\ge 9$ the structure of $\Sigma$ is more or less unknown (even rectifiability is unclear).The singular set
could be a fractal set and will usually have components of varying Hausdorff-dimension $\le n-7$ in $M^{n+1}$.\\

Now we can describe the basic scheme of this paper (and  also of [CL]) as a composition of constructions which we call \emph{cone reducible functors}\footnote{In category theory there are
notions of continuous, asymptotic or tangential functors but with different limit concepts in mind.}: assigning some objects (in our case sets or functions) to $H$ in a way
compatible with cone reductions.\\
 In detail:  Consider $H$ and any of its tangent cones $C$ around some point $p \in \Sigma$ and some construction "K" that assigns some object $K(F)$ to any
 minimal hypersurface $F$ and we assume there is some topology on the space of these objects (for now called "K-topology"), for instance, distinguished functions equipped with the $C^k$-topology (on spaces which can be
 identified via some canonical almost isometric diffeomorphism).\\
 Now we have local flat norm convergence of $\tau_i \cdot H$ to $C$ around $p$ for some sequence $\tau_i \ra \infty$ and we call "K"  \emph{cone reducibly functorial}
 provided that \emph{commutativity} of the following \emph{diagram} holds for any tangent cone $C$ and any such sequence
\begin{center} $
\begin{CD}\tau_i \cdot H@>K>>K(\tau_i \cdot H)\\@VVflat \: normV@VVK-topologyV\\C@>K>>K(C)\end{CD} $\end{center}
which means that the asymptotic behaviour of $K(H)$ near $p \in H$ can be understood from the limit case on $C$. \\
A main result in this paper is that the assignment of some \emph{distinguished} positive eigenfunctions of the conformal Laplacian (note that in the noncompact case there can be many positive
eigenfunctions)  and the involved distance concepts are instances of cone reducible functors. \\

This functoriality will show that conformal Laplacians and their geometric impact via conformal deformations using  first eigenfunctions on singular minimal hypersurfaces can be analyzed
matching naturally with cone reductions leading to a sharp picture of this operator and its eigenfunctions near the singular set. \\ One major application is that for $H \subset (M,g)$ with
$scal_M > 0$ we can find $scal > 0$-metrics on $H$ well-controlled near $\Sigma$ and amenable to stratified surgeries as developed in [CL] which provide a lossless method to eliminate
the singularities. Thus we can incorporate regular and \emph{singular} minimal hypersurfaces on an equal basis
as tools in scalar curvature geometry. \\

Now turning to some technical details, we first point out that the naive strategy to consider the standard conformal Laplacian $L = \triangle + \frac{n-2}{4 (n-1)} \cdot scal_H$ on singular
spaces just as in the smooth compact case, does not lead to a satisfactory theory. Namely, scaling this operator around a singular point does not allow to transfer scalar curvature
information on $H$ to tangent cones backward-and-forward.\\
A way to solve this issue is to redistribute the scalar curvature or equivalently to consider a weighted conformal Laplacian that takes care of the scaling effects. However, to find a suitable
weight we cannot use the usual distance notion (to measure distances between regular points and the singular set $\Sigma$)  since it is not compatible with cone reductions. Instead we will we
introduce another device motivated from the observation that under degenerations from smooth to
singular minimal hypersurfaces one observes  that $| A | $ becomes a measure for a distance to $\Sigma$ which naturally translates to a distance function to the singular set $\sigma$ of any tangent cone.\\
This leads us to consider $L_H = - | A | ^{-2} \cdot ( \triangle + \frac{n-2}{4 (n-1)} \cdot scal_H )$. For the moment let us ignore the fact that $ | A |^{-1}(0) $ may not be empty. Then, for
smooth $H$, $L_H$ will keep the same information as $L$ and its first eigenfunction leads to a conformal deformation on $H$ whose scalar curvature has the same sign as for $L$ and
moreover if $H_i$ is a sequence of smooth hypersurfaces degenerating to some singular $H$ we observe a natural transition to some eigenfunction of $L_H$ of the limit $H$. And in this
singular case $L_H$ will have the versatile feature that its first eigenvalue and eigenfunction carries over (in a way clarified at  length below) to the operators on its tangent cones.  A look at
the eigenvalue equation
\[(\ast) \;\;\;  - \triangle u + \frac{n-2}{4 (n-1)} \cdot scal_H  \cdot  u = \lambda \cdot  | A | ^2  \cdot u  \;  \;  \;  \mbox{ on }H \setminus \Sigma\]
already reveals a key point: namely the scaling invariance of the eigenvalues. Therefore (and motivated from the cone case) we will call this weighted  conformal Laplacian
 $L_H$ also the \emph{scaling invariant conformal Laplacian}.  \\

Now we must have a look at the problem that usually   $ | A |^{-1}(0) \neq \emptyset$. We may assume that $| A |^{-1}(0) \varsubsetneq H$ (otherwise also get $ | A | \equiv 0$ on its
tangent cones and hence $H$ is smooth). Moreover we can assume that it is a nontrivial set of measure zero since we could slightly $C^k$-perturb $(M,g)$ to turn $H$ and thus $A$ into
analytic objects and the tangent cones are (for the same reason) analytic anyway.  After this harmless reduction one has to handle a smaller but now rather persistent set $ | A |^{-1}(0) \neq
\emptyset$. However the upshot is that the problems this  causes can be resolved by approximation
methods developed and explained later on. \\

We consider the equation $(\ast)$ on an area minimizing hypersurface $H$ within  $(M,g)$ with $scal_M > 0$ and on tangent cones $C_p \subset \R^n$ and get\\

{\bf Theorem 1} \quad\emph{Up to multiples there is a unique positive eigenfunction $f_H$ for some eigenvalue $\lambda_H > 1/10$ on $H$. On each \textbf{regular} tangent cone there are
two linear independent positive eigenfunctions for this eigenvalue.\\}

($\lambda_H$ is characterized as the \emph{lim inf} of first eigenvalues of Dirichlet problems on regular domains $\subset \! \subset H \setminus \Sigma$.)\\

 Calling  a function $f_H$ that solves the equation $(\ast)$  an eigenfunction is not quite correct but admissible for our purposes:  $f_H$ is not a first eigenfunction as in the  case of a
 smooth closed $H$. Actually (extending the Martin theory for regular domains $\subset \R^n$,  cf. [P], sec. 4.) we
will find positive functions solving $(\ast)$ for any  $\lambda < \lambda_H$ on $H \setminus \Sigma$ (these $ \lambda$ are called subcritical) and we will use this extensively. Since it is
easily seen that there is no eigenfunction for $ \lambda > \lambda_H$,  $\lambda_H$ is  also called generalized principal eigenvalue of $L_H$.\\

The existence of positive functions solving $(\ast)$ for  $\lambda < \lambda_H$ is also the reason why the case of cones with higher dimensional singular set looks different:\\

{\bf Theorem 1'} \quad\emph{On each \textbf{non-regular} tangent cone $C$ there are infinitely many linear independent positive eigenfunctions for the eigenvalue $\lambda_H$.\\
However there is still a distinguished positive eigenfunction $f_C$ obtained as a limit of first eigenfunctions for Dirichlet problems on regular domains  in $C$.\\}

Specifically, the regular domains in $C$ of Theorem 1'  will be sets of the form $|A|^{-1}([0,a]) $. The smoothness of $\p (|A|^{-1}([0,a]) = |A|^{-1}(\{a\})$ follows from the cone property of
$C$ and Sard's lemma. However if $|A|^{-1}(\{0\}) \neq \emptyset$ it reaches the singular set $\sigma \subset C$ and thus we will actually use some
(again functorial) averaged version  of $|A|$ which is fine enough to be able to \emph{assume} that $|A|^{-1}(\{0\}) = \emptyset$ on $H$ and on its tangent cones (cf. sec. 5 below for some more background).  \\
The second part  of Theorem 1'  fits seamlessly with Theorem 1:  when a non-regular cone $C$ is the limit of a flat-norm converging sequence of regular cones $C_i$ we observe that
$f_C$ can be represented as a limit of a sequence of (also distinguished) positive eigenfunctions on  $C_i$.\\

Now focussing on the limiting behavior of $f_H$ we scale $H$ around a point $p \in \Sigma$ and via Allard regularity there are eventually arbitrarily large compact regular regions on any
tangent cone  which approximate corresponding parts on $H$ in $C^k$-topology. Thus, since our eigenvalue is scaling invariant we consider $f_H$ as a sequence of solutions of $(\ast)$ on
$C_p$
 on growing portions of $C_p$. Imposing some local $L^2$-normalizing this produces (via Harnack inequalities) a positive limit
solution on $C_p$. A priori this \emph{induced solution} on $C_p$ need not to be well-defined.  However this is actually the case: this process selects the one  \emph{minimal} towards the
singular set of the cone. Formally, since we have $scal_H \approx - | A | ^2 $ close to $\Sigma$ and
this minimality "towards" $\Sigma$ is critical only in an arbitrarily small neighborhood of $\Sigma$, we are led to the following definition where we allow also \emph{non-compact} minimal hypersurfaces (like cones).\\

{\bf Definition} \quad \emph{We call a smooth solution $\wp > 0$  of the equation $\triangle \varphi +(\frac{n-2}{4 (n-1)} + \lambda^0) \cdot |A|^2 \cdot \varphi = 0$ a \textbf{Perron
solution} if there is some neighborhood $W$ of $\Sigma$
 with $\p W \cap H \setminus \Sigma$ smooth,  such that $\wp$ is the smallest positive solution of
 $\triangle \varphi +(\frac{n-2}{4 (n-1)} + \lambda^0) \cdot |A|^2 \cdot \varphi = 0, \varphi|_{\p  W } \equiv \wp$ on $\overline{W}$.}\\

This  concept is related to that of \emph{minimal harmonic functions} used for Martin boundaries (cf. [Do] and [P]). However there is an important  difference since our notion of minimality is
 adapted for singularities: as will be shown later it selects a unique element in the Martin boundary which is minimal towards the singular set. Understanding the whole Martin boundary will show that the Perron
solutions on minimal hypersurfaces  can be understood inductively (see Theorem 3 below).\\
 It is notable that the concept of Perron solutions is not just descriptive, used in junction with boundary Harnack inequalities (Carleson inequalities) it becomes a tool to analyze the
Martin boundary in places which had previously been accessible only via probabilistic methods. \\

The name (and the definition) will be justified later by some Perron type (re)construction of solutions on $H$ and its tangent cones. This will also show that if  $\wp > 0$ is minimal with
respect to such a neighborhood $W$ then it is
also minimal with respect to any smaller neighborhood.\\
 Note an important detail: when checking the minimality on $W$ the competing potentially smaller functions need not to be defined/extendable
outside of $\overline{W}$. Thus uniqueness and the Perron property are independent conditions but we will prove\\

{\bf Theorem 2} \quad\emph{$f_H$ and $f_C$ are  the (up to multiples) uniquely determined solutions with Perron property.\\}

 Henceforth $f_H$ and $f_C$ are labelled $\wp_H$ resp. $\wp_C$.\\

Since Theorem 1' shows that the space of solutions can be rather large the question is how these solutions $f_H$ and $f_C$ relate.  This is answered by the following main result of this paper
which says that assigning the Perron solution to a minimal hypersurface is a \emph{cone reducible functor}.\\

{\bf Theorem 3} \quad\emph{$\wp_H$ induces exclusively the Perron solution $\wp_{C_p}$.\\}

In order to derive uniform  estimates (on the space of minimal cones) for the growth and other properties of $\wp_C$ (and to deduce subsequently estimates for $\wp_H$) we use the flat
norm compactness of the set of all minimal cones
and the following (nontrivial) consequence of the theorems above\\

{\bf Corollary 1} \quad \emph{For a flat norm converging sequence of minimal cones $C_i \ra C_\infty$ we have $C^k$-compact convergence $\wp_{C_i} \ra \wp_{C_\infty}$ on smooth
domains (identified
via Allard regularity).}\\

$ \wp_H$ is a limit of a sequence of solutions $u_k$ of Dirichlet eigenvalue problems on  $|A|^{-1}([0,k])  \subset  H  \setminus \Sigma$ for $k \ra \infty$. This leads us  to relate the
uniqueness of $f_H = \wp_H$ and the Martin boundary (at infinity) for $L_C$ of $|A|^{-1}([0,a]) \subset C$ within a given tangent cone.
Actually,  Theorem 3  uses that this Martin boundary is a single point: \\

{\bf Theorem 3'} \quad\emph{There is precisely one positive solution (up to multiples) for the following problem on $|A|^{-1}([0,a]) \subset C$:
 \[\triangle \varphi +\left(\frac{n-2}{4 (n-1)} + \lambda^0 \right)\cdot |A|^2 \cdot \varphi = 0, \; \varphi|_{|A|^{-1}(\{a\}) } \equiv 0.\]}

These theorems allow us to understand the behavior of $\wp_H$ via the shape of $\wp_{C_p}$ which can be analyzed by some induction scheme.
Writing points in $C_p$ in polar coordinates, i.e. $(\omega,r) \in C_p$ where $r $ is the distance to the tip and $\omega$ a point in $\p B_1(0) \cap C_p$\\

{\bf Theorem 4} \quad\emph{$\wp_{C_p}$ admits a separation of variables: $\wp_{C_p} = c(\omega) \cdot r^\alpha$ for some positive function $c(\omega) > 0$ solving
\[(CW) \;\; \left(\alpha^2 + (n-2) \alpha \right) \cdot c(\omega) +
\left( \triangle_S + \left( \frac{n-2}{4 (n-1)}+\lambda \right) \cdot a(\omega)^2 \right) c(\omega) = 0 \] $-\frac{n-2}{2}  < \theta_1(n) < \theta_2(n) < 0$ such that
$\alpha \in (\theta_1(n),\theta_2(n))$ and $\alpha = \alpha_p$ is uniquely determined for every $p \in \Sigma$}.\\

Near the singularities of $C_p \cap \p B_1(0)$ we observe that after scaling $(CW)$ becomes again an equation of the form $(\ast)$ however with dimensions shifted (which can be handled
like the
$(n-1)$-dimensional form of $(\ast)$) and we get from Theorem 1 and 2: $c(\omega)$ is again the unique positive solution and has the Perron property.\\

Beside the fact that Theorem 4 and Corollary 1 enter in the proof of the previous results via induction we note some other important consequences for  $\wp_H^{4/n-2} \cdot g_H$ and
$\wp_{C_p}^{4/n-2} \cdot g_{C_p}$. Since the Perron solution has a growth near
$\Sigma$ as $(r^\alpha)^{2/n-2}$ for $- \frac{n-2}{2} < \alpha$ the length function has an integrable singularity in  $0 \in \R^+$:\\

{\bf Corollary 2} \quad\emph{The diameter of $(H,\wp_H^{4/n-2} \cdot g_H)$ is finite. } \\

(Actually there is some uniform control discussed later on.)\\

The applications to scalar curvature geometry are combinations of Theorem 3 and of the following consequences of Theorem 4 saying that after applying the conformal deformation
$\wp_H^{4/n-2} \cdot g_H$ the
geometry near a point in $\Sigma$ looks like a cone with $scal \ge 0$.\\

First note that any space of the form $N \times \R^{\ge 0}$ with a warped product metric $ r^2 \cdot g_N + g_{\R}$  where $(N, g_N)$ is an arbitrary Riemannian manifold,
is an (abstract) cone (i.e. scaling invariant around $0$) and vice versa such a scaling invariant space can be written as a warped product.\\

{\bf Corollary 3} \quad\emph{\begin{itemize} \item Any cone $(C_p,g_C)$ equipped with the metric $\wp_{C_p}^{4/n-2} \cdot g_C$ is again a cone and it has $scal \ge 0$. \item The metric
$\wp_{C_p}^{4/n-2} \cdot g_C$ is conformal to another cone metric $\overline \wp_{C_p}^{4/n-2} \cdot g_C$ with $scal_{ \overline \wp_{C_p}^{4/n-2} \cdot g_C}(\omega,\rho) \ge
\iota_H/\rho^2$ for some $ \iota_H > 0$ independent of $C_p$. \end{itemize}}

($\rho$ is the distance to the tip of the cone with respect to $\overline \wp_{C_p}^{4/n-2} \cdot g_C$)\\

There is also a corresponding deformation from  $\wp_H^{4/n-2} \cdot g_H$ to $\overline \wp_H^{4/n-2} \cdot g_H$  done in the same functorial way  that commutes with the transitions
to tangent cones. And clearly this needs the scaling invariant version of the conformal Laplacian. The advantage is  when we can  deform the metric $\wp_H^{4/n-2} \cdot g_H$ additionally in
a way that allows us to perform some kind of surgery or some other regularization process close to $\Sigma$ we can now zoom into as deep as we want (to gain local simplifications of the
geometry) without loosing the local
positive lower bound for $scal$ needed to compensate for additional deformations carried out in that region during any sort of regularization.  \\

\vspace{1.7cm}
 \setcounter{section}{2}
\renewcommand{\thesubsection}{\thesection}
\subsection{Area minimizing cones and reduction techniques} \label{area-min-cones}

\bigskip

The only a priori information concerning the singular set $\Sigma \subset H^n$ we use is the compactness and the Hausdorff-dimension which is $\le n-7$. But we have a structural aid provided by tangent cones (cf.~\cite{Gi}, \cite{Si}). These (locally area minimizing) minimal cones in $\R^n$ are a generalization of the tangent plane at singular points: after some scaling one may consider $H$ as locally (say around $p \in \Sigma$) embedded in $\R^n$, and after further scalings by an increasing sequence of factors $\tau_m \ra +\infty$ there is a minimal cone $C_p$ which approximates $\tau_m \cdot H$ on any given compact set in $\R^n$ in a certain way described below.\\

The usage of tangent cones in the literature is fairly limited since each singular point in $\Sigma \subset H$ will usually have infinitely many tangent cones, the set of tangent cones varies discontinuously along $\Sigma$ and the approximation of $H$ by these cones is {\it not} uniform in $\Sigma$. \\

Nevertheless, because we will avoid to come \emph{too} close to $\Sigma$, we will be able to set up a scheme to derive many properties of $H$ near $\Sigma$ from corresponding information on cones. For certain properties this even allows us to gain \emph{uniform} control by using the precompactness of the space of tangent cones. Finally, and most importantly, the approximation by tangent cones will allow us to carry out certain local operations on cones (serving as models), and then transplant them to $H$.\\

We start on an abstract level with a composition of several classical facts due to De Giorgi, Allard and others (cf.~\cite{DG}, \cite{A1}, \cite{Gi}, and \cite{Si}). \\

\begin{proposition} \label{flat-norm-approx}
\quad  Let $H^n \subset M^{n+1}$ be an area minimizing hypersurface and $\tau_m \to +\infty$ a sequence of positive real numbers.\\
Then, for every $p \in \Sigma$ we find a subsequence $\tau_{m_k}$ and an area minimizing cone $C_p \subset \R^{n+1}$ such that for any given open $U \subset \R^{n+1}$ with compact closure the {\it flat norm} $d_U$ (cf.~\cite{Si}, Ch.~31) which (roughly speaking) measures the volume between two sets in U  converges to zero:
$$
d_U (\tau_{m_k} \cdot H, C_p) \to 0.
$$
Moreover, if $\overline U$ contains only smooth points of $C_p$, this convergence implies {\it compact $C^l$-convergence}, for any $l \ge 0$.
\end{proposition}

\begin{remark} \label{flat-norm-approx-remark}
\quad Using normal coordinates $\tau_{m_k} \cdot H \subset \tau_{m_k} \cdot M$ can locally near $p$ be considered as a subset of $\R^{n+1}$ (for $k \to +\infty$ the deviation vanishes). The $C^l$-convergence statement can be obtained by combining Allard regularity with elliptic regularity. It can be formulated more precisely as follows: let V be an open subset of $C_p$ with focal distance $\iota > 0$, whose compact closure contains only regular points. Consider the $exp_\nu$-image $U_\varepsilon$ of normal vectors of length $\le \varepsilon \le \iota/2$ in the normal bundle $\nu|_V$ over $V$. Then for large $k$ the set $U_\varepsilon \cap \tau_{m_k} \cdot H$ is a $C^l$-graph (= $C^l$-section of the normal bundle) over $V$, and converges compactly to $V$ (= zero section) in $C^l$-topology.
\end{remark}

The cone reduction argument we are looking for cannot be based on particular properties of a special cone, but becomes valid only if we can manifest such properties for the class of all singular cones simultaneously. One of the ingredients will therefore be the following two results.

\begin{lemma} \label{cone-compactness}
\quad  The set $\mathcal{C}_{n}$ of embedded area minimizing $n$-cones (around $0$) in $\R^{n+1}$ is compact in the flat norm topology.
\end{lemma}

{\bf Proof} \quad This can be derived from the compactness theorem for
integral currents (see e.g.~37.2 in \cite{Si}), and the fact that minimality and the cone shape survive under flat norm convergence. \qed

In particular, the set  ${\cal T}_H$ of singular tangent cones of $H$ (with center set to 0) has the \emph{compact} closure $\overline {\cal T}_H \subset\mathcal{C}_{n}$. $\overline {\cal T}_H$ will usually contain cones which do {\bf not} appear as tangent cones of $H$. Actually considering such extensions deliberatively will be an essential tool for many arguments. We state a simple but crucial compactness result in this direction:

\begin{corollary} \label{singularity-criterion}
\quad There is a constant $d_n >0$ such that
$$d_{B_1(0) \setminus B_{1/2}(0)} (C, \R^{n}) < d_n \;\;\; \text{if and only if} \;\;\; C \;\text{is non-singular.}$$
Therefore the set of singular cones  $\mathcal{SC}_n \subset\mathcal{C}_n$ is closed (and hence compact). Consequently $\overline {\cal T}_H \subset \mathcal{SC}_n$.
\end{corollary}

{\bf Proof} \quad Let $C_i$ be a sequence in $\mathcal{SC}_n$ with $d_{B_1(0) \setminus B_{1/2}(0)} (C_i, \R^{n}) \to 0$. Then by the cone property $d_{B_1(0)}(C_i,\R^{n}) \to 0$, and Allard regularity implies that for large $i$ every $C_i$ is non-singular. \qed

Next we will sharpen the usual picture of cone approximation: For decreasing radius $\eta \ra 0$ $(\eta^{-2} \cdot H) \cap B_2(p) \setminus B_1 (p)$ is not just sometimes approximated by
a cone but a slightly closer look already unveils an instructive view: choose a
\emph{finite} covering $\{ B_\delta (c_i) \}$ of the compact set of singular cones $\cal{CS}$ by \emph{flat norm} balls of radius $\delta$.\\
1. For any $\delta
> 0$ we find that starting from some $\eta_\delta
> 0$ such that $(\eta^{-2} \cdot H) \cap B_2(p) \setminus B_1 (p)$ is $\delta > 0$ - close in flat norm to some (non uniquely determined) tangent cone
$C^\eta_p$.\\
 2. Considering this assignment as a discrete valued map $\eta \mapsto \{ B_\delta (c_i) \}$ we observe a large scale \emph{fading} or \emph{freezing} property:
 after scaling $\eta$ to 1 the frequency of oscillation within
the balls of this finite covering will decay uniformly to zero for  $\eta \ra 0$ and (also after scaling) the size of the well-approximated part
of \emph{any} of these cones increases (i.e. considering a sequence of approximating regions (identified via scaling) we get a compact exhaustion of any tangent cone).\\

This is just an interpretation of the following \\

\begin{lemma} \label{approx-lemma}
\quad For any $\delta > 0$ and any $R \gg 1 \gg r >0$ we can find a small $\eta_{\delta , R , r} > 0$ such that for \textbf{every} $\eta \in (0, \eta_{\delta , R , r})$ and any tangent cone
$C_p^\eta$ of $H$ at $p$:
$$(\eta^{-2} \cdot H) \cap (B_R(p) \setminus B_r (p)) \text{ is $\delta$-close in flat norm to } C^\eta_p\cap(B_R(0) \setminus B_r (0)).$$
\end{lemma}

Note that $\eta_{\delta , R , r}$ depends on $p$ in a discontinuous way.\\

The \textbf{proof} is standard: if there were a sequence of $\eta_i \ra 0$ and a $\delta_0 > 0$ such that $\eta_i^{-2} \cdot H \cap B_R(p) \setminus B_r(p)$ is not $\delta_0$-close to any tangent cone, there is still a subsequence that gets arbitrarily close to some tangent cone.\qed

Notice that Allard regularity combined with elliptic regularity provides us with the refined version for $C^k$-topology: Suppose $\sigma^\eta_p$ denotes the singular set of some tangent cone $C^\eta_p$, and $V_{a}(\sigma^\eta_p)$ the union of the {\it subcone} in $C^\eta_p$ over a sufficiently small neighbourhood $U_{a}(\sigma^\eta_p\cap \p B_1(0)) \subset  (C^\eta_p\cap\p B_1(0))$ with $B_a(0)\subset C^\eta_p$. Then \ref{approx-lemma} implies together with Remark \ref{flat-norm-approx-remark}:

\begin{corollary}\label{aprr}
\quad For any $\delta > 0$ and any triple $R \gg 1 \gg r \gg a >0$ we can find a small $\eta_{\delta , R , r, a}> 0$ such that for {\it every} $\eta \in (0, \eta_{\delta , R , r, a})$ the corresponding part of $\eta^{-2} \cdot H$ can be written as the graph of a function $g_\eta$ over $C^\eta_p\cap B_R(0)\setminus (B_r (p)\cup V_{a}(\sigma^\eta_p))$ such that $ |g_\eta|_{C^k} < \delta$. \qed
\end{corollary}

This suggests an important relation "$\Sigma \prec \sigma$" between the singular sets $\Sigma \subset H$ and $\sigma \subset C_p$ : asymptotically the singular set of the tangent cones is  "larger" than the germ of the singular set around $p \in H$. For instance, $\Sigma$ may contain scattered points, or there might be smooth but highly curved regions near $\Sigma$ which may cause the appearance of rays in $\sigma$ . On the other hand, however, the complexity of any $\sigma$ is reduced by one dimension (since $\sigma$ is also a cone).\\

Later on this fact will play a crucial role: The conformal deformations close to $\Sigma$ that we use to form a barrier around $\Sigma$ will first be prepared on tangent cones (instead of $H$) (see Section \ref{local-barriers}), and then transplanted to the regular regions of sufficiently well-approximated balls in $H$ (see Section \ref{global-barriers}). In Section \ref{besicovich-coverings} we will show in detail how these balls can be obtained. Thus, the constructed barriers will hide not only $\Sigma$, but rather a whole neighbourhood of $\Sigma$ that is induced, and in some sense stratified, by very small neighbourhoods of the $\sigma$s of the approximating tangent cones. In the following we will refer to this stratification of $\Sigma$ as the {\it local enhancement} of $\Sigma$.\\

We now describe the basic procedure we use to mediate between $H$ and the realm of singular cones and how to proceed from there.\\

In order to prove a local result on $H$ which is known to be true for cones we frequently argue by contradiction. Assume there is a sequence of points $x_n \in H \setminus \Sigma$, $d_M(x_n, \Sigma) = \ve_n \to 0$ (Note that we will consider the intrinsic distance later on) and around $x_n$ a certain expected
geometric (or more general analytic) property fails to hold on $B_{\alpha \cdot \ve_n}(x_n)$, $\alpha \ll 1$. In addition, the property in question should satisfy a compactness property: e.g.~elliptic compactness (and Arzela-Ascoli) when we consider eigenfunctions, Gromov compactness (plus Allard regularity) when we consider the second fundamental form as a curvature quantity. \\
We will then argue as follows: There is a $p \in \Sigma$, being limit of a subsequence of $x_n$, and $\rho_n = d_M(x_n,p) \ge \ve_n$ will also converge to zero. There are two cases
\begin{enumerate}
\item $\ve_n / \rho_n > const. > 0$: in this case the $x_n$ run into a well-approximated zone of a tangent cone in $p$,
\item $\ve_n / \rho_n \ra 0$: here we still get a cone approximation, but the cone need not appear as a tangent cone at any point of $\Sigma$.
\end{enumerate}

In case (i), after scaling $H$ and $M$ by $\ve_j^{-2}$ (so that $d_M(x_j,p)$ is normalized to 1 up to bounded multiple), there is still a subsequence of $x_{j_k}$ converging (in this scaled picture) to a point $q \in \p B_1(0)\cap C_p$ where $C_p$ is a tangent cone at $p$.\\

In case (ii), we can argue as follows: Take a point $p_j \in \Sigma$ with $d_M(x_j,p_j) = d_M(x_j, \Sigma) =  \ve_n$  and scale each intersection $H \cap B_{ \rho_j}(p_j)$ by $\rho_j^{-2}$. This can be considered a sequence of area minimizing surfaces $T_j$ in $B_1(0) \subset \R^{n+1}$, and we may assume it converges in flat norm to an area minimizer $T_\infty$ in $B_1(0) \subset \R^{n+1}$. Now, a subsequence of $(\ve_j / \rho_j)^{-2}$-scaled copies of $T_\infty$ converges in flat norm to a minimal cone $C_\infty$ (which need not be a tangent cone of $H$). Thus, by a diagonal sequence argument we may assume that  $H \cap B_{ \rho_j}(p_j)$ scaled by $\ve_j^{-2}$ converges in flat norm to $C_\infty$ and that $x_{j}$ converges (in this scaled picture) to a point $q \in \p B_1(0) \cap C_\infty$.\\

For  convenience we will use $C_\ast$ as a common notation for $C_p$ resp.~$C_\infty$ when both cases can show up. The second case will be called an \emph{abstract cone reduction}. In either case the limiting cone is smooth outside a \emph{codim 7} singular set $\sigma$, hence the flat norm convergence gives rise to compact $C^l$-convergence outside $\sigma$.\\

Now, the cone reduction strategy proceeds as follows: In certain cases an \emph{a posteriori} argument shows that $q$ is a \emph{regular point} in  $C_\ast$. In some other cases we use that after scaling around $q$, $C_\ast$ can be approximated by a tangent cone which is a product $\R \times \hat{C}^{n}$, where $\hat{C}^{n} \subset \R^{n}$ is again a minimal cone and argue inductively. Then we may use the compactness result for the geometric/analytic estimate or property under consideration and the fact that $\ve_n^{-2} \cdot H$ converges to $C_\ast$ to conclude the estimate/property continues to fail on $B_\alpha(q) \subset C_\ast$. Therefore we are done if we know that, in fact, the corresponding property does hold on $B_\alpha(q)$. \\

Direct arguments (and hence sharper estimates) often fail since this would usually require uniform approximation by tangent cones. \\

The cone structure actually provides us with two tools: firstly, the cone direction which blows up to give a local product structure with a minimal hypersurface $G^{n-1}$ as a construction aid on its own, and secondly, the properties of $G^{n-1}$ that can be used as hypothesis for the next step of the induction. \\

\begin{remark} \label{cone-sphere-intersection}
\quad At this point it is important to mention that $G^{n-1} = \p B_1(0) \cap C$ is minimal but neither area minimizing nor stable (since $\p B_1(0)$ has $\Ric >0$). Nevertheless we can carry over those results valid for area minimizers that allow us to make the induction work: The crucial property of $G^{n-1}$ in this setting is that the \emph{cone over} $G^{n-1}$ is area minimizing, and hence all its tangent cones are. Outside 0 the tangent cones have a product structure isometric to $\R \times \tilde C_q^{n-1}$ where $\tilde C_q^{n-1}$ is again area minimizing. But these cones $\tilde C_q^{n-1}$ are precisely the tangent cones of $G^{n-1}$. This, together with the local product structure of $C$ as a cone over $G^{n-1}$ will allow us to handle $G^{n-1}$ is our scheme just like an actual area minimizer, e.g.~the singular set of $G^{n-1}$ has the same properties (e.g. ${\rm codim} \ge 7$, compactness) as that of area minimizers. In addition, the argument for distinguishing the two cases for $C_\ast$ survives. This would not be the case for general minimal surfaces. \\

Also, there will be no accumulating problem during the induction process. In the next step we consider the tangent cones $\tilde C_q^{n-1}$ of $G^{n-1}$, consider $\p B_1^{n-1}(0) \subset \R^{n-1}$ and $G^{n-2}$ etc.~until we obtain isolated point singularities. \\

A good way of thinking of this part of the strategy is as a more complex version of the classical cone reduction used to determine the codimension of $\Sigma$ -- just with additional data on the hypersurfaces inducing corresponding data on the lower dimensional objects.
\end{remark}

As a sample of this rather abstract scheme we consider the \emph{intrinsic} distance function on $H$. The extrinsic distance between points $x \in H \setminus \Sigma$ and (points in)
the compact set $\Sigma \subset M$ measured within the ambient manifold $d_M(x,p)$ resp. $d_M(x, \Sigma)$ is not suitable for our purposes: We use the intrinsic metric on $H$ to study e.g.
eigenfunctions of the conformal Laplacian. Moreover, for our argument we will conformally deform the induced metric on $H$
 and, thereafter, we want to understand the the new geometry near $\Sigma$. Yet, at that stage the embedding has lost its meaning. Hence, we have to work with the intrinsic distance function $d_H$ on $H$.\\
Since $H$ may develop additional bumps and even new topology when approaching $\Sigma$ (reflected by thin regions with large $|A|$) one realizes that it is not at all clear that $d_H(p,x) < +\infty$ for any points $x \in H\setminus \Sigma$, $p \in \Sigma$ in the same connected component of $H$.\\
However, using the fact that $H$ is an area minimizer we will prove below that close to $p$ there is a network of pieces of rays which link $x$ to $p$ in finite distance. We will base this argument on a cone reduction. \\

\begin{corollary} \label{intrinsic-metric}
\begin{enumerate}
\item Let $p \in \Sigma$. Assume that $B_r(p)\cap H$ is connected for all sufficiently small $r>0$. Then $(B_r(p) \cap H) \setminus \Sigma$ is also connected, and $d_H(p,x) \le c \cdot r$ for all $x \in B_r(p)\setminus \Sigma$.
\item If $H$ is connected, then $H \setminus \Sigma$ is connected, too. Moreover, its intrinsic diameter is finite.
\end{enumerate}
\end{corollary}

(The ray-network argument we have chosen here to prove the corollary extends directly to situations where we conformally deform $H$, and recover the new intrinsic distances by considering the induced geometries on the tangent cones, e.g.~in Section \ref{}.)\\

{\bf Proof} \quad If $\Sigma$ is a finite set (i.e.~all tangent cones are regular with singularities only in 0) then the tangent cones are connected because of the maximum principle. Due to codimension $\ge 2$, removing $\{0\}$ keeps the complement connected. For $r>0$ small enough we can assume that the set $(B_2(0) \setminus B_{1/2}(0)) \cap C_p$ is $C^k$-close to $(r^{-2}\cdot H)\cap  B_2(p) \setminus B_{1/2}(p)$ for a suitable tangent cone $C_p$ at $p\in\Sigma$. Thus $(B_r(p) \cap H) \setminus \Sigma$ can be written as a union of connected sets (namely rescaled versions of $(B_2(p) \setminus B_{1/2}(p)) \cap H$). Hence, it is connected. The claims concerning intrinsic distances are easily checked in this case. \\

Now proceed with the case where the tangent cones may also contain singularities other than 0. First of all, we claim that for any given $p \in \Sigma$ there is a $k_p>0$ such that for any $r >0$ small enough, $(B_{2r}(p) \setminus B_{r/4}(p)) \cap H$ contains an open connected subset $V_r$ of approximate subcone-shape with the following properties:
$${\rm diam}_{V_r} V_r \le k_p \cdot r, \quad \frac{{\rm vol} (V_r \cap \p B_{2r}(p))}{{\rm vol} (H \cap \p B_{2r}(p))} > \frac{3}{4}, \quad \frac{{\rm vol} (V_r \cap \p B_{r/4}(p))}{{\rm vol} (H \cap \p B_{r/4}(p))} > \frac{3}{4},$$
where ${\rm diam}_{V_r} V_r$ is the intrinsic diameter of $V_r$.\\

To prove this claim, assume that each $k>0$ has a $r_k>0$ such that $B_{2r_k}(p) \setminus B_{r_k/4}(p)$ does not contain such a $V_{r_k}$ satisfying the conditions for $k$ in the role of $k_p$. Since $\Sigma$ has higher codimension, we can assume that the condition on ${\rm diam}_{V_{r_k}} V_{r_k}$ is the one that fails. Passing to a subsequence, we may assume that $r_k^{-2}\cdot(B_{2r_k}(p)\setminus B_{r_k/4}(p))$ converges to $C_p\cap(B_2(0)\setminus B_{1/4}(0))$ for some tangent cone $C_p$. Since $B_{r_k}(p) \cap H$ is connected by hypothesis, $\p B_1(0) \cap
C_p$ is connected, too. So, arguing inductively, we may assume the claim to be true in codimension 1 (cf.~ Remark \ref{cone-sphere-intersection}). Hence, let us assume that $\p B_1(0) \cap C_p \setminus \sigma$ is connected, where $\sigma$ is the singular set of $C_p$, that the intrinsic diameter of $\p B_1(0) \cap C_p$ is finite, and that there is an open connected $\tilde{W} \subset (\p B_1(0) \cap C_p) \setminus \sigma$ with $\frac{{\rm vol} (\tilde{W})}{{\rm vol} (C_p \cap \p B_1(0))} > \frac{4}{5}$, ${\rm diam}_{\tilde{W}} \tilde{W} < +\infty$. However, defining $\tilde V_{r_k}=(B_{2r_k}(0) \setminus B_{r_k/4}(0)) \cap (\mbox{subcone over } \tilde{W} \text{ in } C_p)$, and using the $C^l$-approximation of the scaled $H$, we obtain a corresponding set in $H$ giving a contradiction.\\

Note that the compactness result for tangent cones allows us to adjust the $V_r$ such that ${\rm diam}_{\tilde{W}_r} \tilde{W}_r$, where $\tilde{W}_r:=\p B_1(0) \cap \tilde{V}_r$ for the corresponding sets $\tilde{V}_r$ in $C_p^r$, is uniformly bounded from above by some constant $b<\infty$. Otherwise there is sequence of tangent cones $C_p^{r_j}$ converging to some $C$ (in $C^l$ on smooth parts) such that of ${\rm diam}_{\tilde{W}_{r_j}} \tilde{W}_{r_j}$ diverges etc.\\

Now let $r>0$ be small, and consider $\cup_{k=0}^\infty V_{r/2^k}$ (note that there are several tangent cones involved). Because the volume fraction belonging to $V_{r/2^k}$ in both $\p B_{r/2{k-1}}(p)$ and $\p B_{r/2{k+2}}(p)$ is larger than $\frac{3}{4}$, and $V_{r/2^k}$ has approximate sub-cone shape, the intersection  $V_{r/2^k} \cap V_{r/2^{k-1}}$ is open and non-compact in $\p B_{r/2^{k-1}}(p)$. Starting at $x$ we now choose a path which follows the (approximate) ray direction in $V_r \cap (B_r(p) \setminus
B_{r/2}(p))$. Then, on $\p B_{r/2}(p)$ one uses ${\rm diam}_{\tilde{W}_{r/2}} \tilde{W}_{r/2} < b$ to run to a point $x_1$, from which one can follow an approximate ray direction within $V_{\rho/2}$ to reach an intersection point with $V_{\rho/4}$ etc. Thus we get a sequence of points $x_k \in H \setminus \Sigma$ with $x_k \to p$ and $d_H(x_k,x_{k+1}) \le (1+b)\cdot 2^{-k} \cdot r$. Thus, $d_H(p,x) \le (1+b)\cdot r$ when $x \in B_r(p)$. The other claims are direct consequences of this construction. \qed

\begin{corollary}
\quad There is are universal bounds $0 < A_1(n) < A_2(n) < \infty$ and $0 < D_1(n) < D_2(n) < \infty$ for the area $A$ and diameter ${\rm diam}$ of $\p B_1(0) \cap C$ for any $C \in SC_n$:
$$A_1(n) < A < A_2(n) \;\;\;\;\;\mbox { and }\;\;\;\; D_1(n) <  {\rm diam} <  D_2(n).$$
\end{corollary}

{\bf Proof} \quad This is a consequence from the compactness of $\mathcal{SC}_n$, and we only indicate the argument for the least obvious claim ${\rm diam}_{\p B_1(0) \cap C} < D_2(n)$. If $C_j$ is a sequence with ${\rm diam} _{\p B_1(0) \cap C_j} \ra \infty$ we may assume it converges in flat norm and $C^k$-compactly to some limit cone $C_\infty$. According to \ref{intrinsic-metric}, however, we have ${\rm diam} _{\p B_1(0) \cap (C_\infty\setminus\sigma_\infty)} =: D < \infty$. Hence, the compact $C^k$-convergence implies that there is a sequence $\ve_j \ra 0$ such that ${\rm diam} _{\p B_1(0)\cap C_j \cap U_{\ve_j}} \ra \infty$ where $U_{\ve_j}$ is the extrinsic $\ve_j$-neighborhood of $\sigma_j$. Rescaling by $\ve^{-2}_j$ as $j\to\infty$, we approximate a minimal hypersurface (cf.~Remark \ref{cone-sphere-intersection}). However, we can assume uniform diameter bound for this minimal hypersurface, by using \ref{intrinsic-metric} together with an inductive cone reduction argument.\qed

\textbf{Remark} \quad In what follows we can therefore assume that $H$ and $H \setminus \Sigma$ \emph{are} connected, since the subsequent arguments will apply to each component.\\

\vspace{1cm}

\setcounter{section}{3}
\renewcommand{\thesubsection}{\thesection}
\subsection{Distance Functions on Minimal Hypersurfaces}

\bigskip

The aim of this section is to find natural distance notions adequate for a fine analysis of/on a singular minimal hypersurface $H$ close to its singular set $\Sigma$. \\

First note that the metric distance of a regular point to $\Sigma$ is \emph{not} helpful  since this distance function does not converge to the distance function to the singular
set in its
tangent cone under scalings. \\
To get on the right track we write a cone  in coordinates $(\omega,r) \in \p B_1(0) \cap C \times \R^{\ge 0} \cong C$ and notice that $|A|(\omega,r) = a(\omega) \cdot r^{-1}$. Thus
$|A|(\omega,r)$ could be used as a distance measure (between regular points and the singularity) except for the case where $a(\omega) = 0$. Before we delve into the nature of
this defect, let us look at the possible advantages of using $|A|$ as a means to measure the distance to the singular set: clearly for $c \ra \infty$, $|A|^{-1}([c,\infty)) \subset H$
shrinks to the singular set $\Sigma \subset H$ and  by the virtue of Allard regularity we see that  for large $c \gg 0$ (some component of) the level sets of $|A|$, $|A|^{-1}(c)
\subset H$ coincide with level sets in $C$ (locally and up to any prescribed
precision).\\
This cone transition of level sets is essential when we use tangent cones as a tool to understand the analysis on $H$ inductively (note this can be iterated until we reach a
regular cone)
and this completely fails if one takes metric distance sets instead.\\
Moreover $|A|$ measures the flatness of the underlying space and this shows that $|A|^{-1}([0,c]) $ has uniformly bounded geometry and from this we get uniform estimates for
elliptic regularity result which survive cone reductions.
Once again the uncontrollable structure of the singular set also diminishes such estimates for metric distance sets. \\

In short, $|A|$ appears to be a natural distance notion near the singular set of minimal hypersurfaces. However, so far we deliberately ignored the fact  that $|A|^{-1}(0) $ may be
nonempty on $H$ or some of its tangent cones. In this case, $|A|^{-1}(c) $ may reach $\Sigma$ which is clearly unwelcome.
Thus our task is to find modifications of $|A|$ which share the advantages of $|A|$ (in particular the cone reducibility of the definitions) but  whose levels behave properly.   \\
Firstly $| A |^{-1}(0) \varsubsetneq H$ (otherwise also get $ | A | \equiv 0$ on its tangent cones and hence $H$ is smooth), next we can assume that $| A |^{-1}(0)$ is  a set of
measure zero since we could slightly $C^k$-perturb $(M,g)$ to turn $H$ and thus $A$ into analytic objects. Tangent cones are (for the same reason) analytic anyway. However
after this
 reduction there is a smaller but persistent set $ | A |^{-1}(0) \neq \emptyset$ on $H$ and/or inductively its tangent cones. \\
Thus we enhance the whole device with some smoothing or, viewed differently, averaging technique utilizing minimal hypersurfaces $N^{n-1}$ within $H^n$ \emph{with obstacles}
equal to level sets of $|A|$ as obstacles. The resulting function ${ }^\approx |A|$ will
be the distance function (towards $\Sigma$) which serves as a natural and cone reducible substitute for the metric distance to $\Sigma$.\\

Thus we briefly digress on \emph{parametric minimal hypersurfaces with obstacles}. In the non-parametric (= graph-type) case  one has global $C^{1,1}$-regularity (cf. [C] for a
reference). But we clearly have to consider parametric hypersurfaces. A priori one has the same type of \emph{codim 7}-singularities as in the case without obstacles while at
least close to the coincidence set
(with the obstacles) one also has $C^{1,1}$-regularity (cf. [M],[KS],[T] and [SZW]) including the usual types of compactness results for free minimal hypersurfaces. \\

 Formally, take two (for now) smooth compact (or complete) and cobordant but not necessarily connected submanifolds $M_1^{m}, M_2^{m}$ and the cobordism
$W^{m+1}$ equipped with some Riemannian metric.
\begin{definition}
 \quad \emph{A (locally) area minimizing current ${\cal{J}}$ in $W^{m+1}$ homologous to $M_1^{m}$ (and thus to $M_2^{m}$) is called an
area minimizer with obstacles $M_1^{m}$ and $ M_2^{m}$.}\\
\end{definition}

In applications one of the two obstacles (say $M_2^{m}$) will usually just be a replacement for a local compactness condition and is never really touched by the support of
${\cal{J}}$  (it may be conceived being placed close to infinity) and thus we will only refer to the effective obstacle as \emph{the} obstacle. More concretely, we will place tiny
neighborhoods $V$  around $\Sigma \subset H$ and consider them as obstacles. To prevent the area minimizer from just collapsing to a point (note that $\p V$ is null-cobordant)
we
will always presume that the area minimizer has to stay in another much larger neighborhood surrounding $V$. \\

We will now observe that $|A|^{-1}(c) $ can serve as an obstacle although it is not complete since the places where completeness fails will be out of reach for area minimizers. To
check this claim we use the auxiliary hybrid $\pounds^2_\ve (x) := \frac {\ve^2}{dist_H(x,\Sigma)^2} + |A|^2(x)$.
 This function is Lipschitz but will not
be smooth in general, but letting the heat flow slightly deform this function gives a smooth approximation (which can be made arbitrarily fine when approaching $\Sigma$) and
henceforth we think of such a fine smooth approximation when we speak of level sets $\pounds_{\ve}^{-1}(d)$ which therefore can (generically) be assumed to be \emph{smooth})
and note the essential
properties that $dist(\pounds_{\ve}^{-1}(d), \Sigma) \ge \ve/d$ and $\ch^{n-1}(\pounds_{\ve}^{-1}(d), \Sigma) \ra 0$ for $d \ra \infty$. ($\ch^{k}$ denotes the $k$-dimensional Hausdorff-measure).\\

Now we want to use $\pounds_{\ve}^{-1}(d)$ as an obstacle. More precisely, we mean the \emph{outermost components} ${}^{out}\pounds_{\ve}^{-1}(d)$ separating $\Sigma$
from the path component of $H \setminus {}^{out}\pounds_{\ve}^{-1}(d)$ which contains (for very large $d$) almost  all of the total volume. Note from [CL],(2.7) that we may
assume that $H$ and $H \setminus  \Sigma$ are connected.\\

\begin{lemma}\label{epsi}
 \quad Let $U \subset H^n$ be any neighborhood of $\Sigma$, and $\ve > 0$ be fixed. Then there is a $d_0
> 0$ such that for almost every $d>d_0$: ${}^{out}\pounds_{\ve}^{-1}(d) \subset U$ and each area minimizer
$h_{d,\ve}^{n-1}$ with obstacle ${}^{out}\pounds_{\ve}^{-1}(d)$ satisfies $h_{d,\ve}^{n-1} \subset H^n \cap U$.\\
\end{lemma}

\textbf{Proof} \quad Of course, since  $\pounds$ is continuous, ${}^{out}\pounds_{\ve}^{-1}(d) \subset U$ for large $d$. On the other hand, we can find for a sequence of
neighborhoods $U(k) \supset U(k+1)$ of $\Sigma$ with smooth boundary and  $\ch^{n-1}(\p U(k)) \le 1/k$. (This is just the co-area formula plus the definition of the Hausdorff
measure.) Choose $d_j$ such that ${}^{out}\pounds_{\ve}^{-1}(d_j) \subset U(j)$ Then, the area minimizing hypersurface $h_j$ with obstacle ${}^{out}\pounds_{\ve}^{-1}(d_j)$
homologous to $\p U(j)$ (both are level sets can be assumed to be smooth) will have $\ch^{n-1}(h_j) \le \ch^{n-1}(\p U(j)) \le 1/j$. Now say there is a $k_0$ such that $h_j \cap \p
U(k_0) \neq \emptyset$. Via compactness we can adapt the situation such that there is common point $p_0 \in h_j \cap \p U(k_0)$ and thus for large $j$ $supp(h_j) \cap
B_{k_0/2}(p_0)$ is a free area minimizer with boundary data on $\p B_{k_0/2}(p_0)$ and $p_0 \in supp(h_j)$. But $B_{k_0/2}(p_0)$ is ball with fixed (bounded) geometry and
thus we had independently of $j$: $\ch^{n-1}(h_j) \ge \ch^{n-1}(supp(h_j) \cap B_{k_0/2}(p_0))> const.>0$.
 \qed

\medskip
 Now we send $\ve \ra 0$: if there is a path or sequence of points in $H \setminus \Sigma$ converging to $\Sigma$ such that $|A|$ along this route converges to a finite
 value, then ${}^{out}\pounds_{\ve}^{-1}(d)$ for $\ve \ra 0$ approaches $\Sigma$ in some points. However $h_{d,\ve}^{n-1}$ does not follow but remains outside a neighborhood $V_d$ of $\Sigma$
 cf. Proposition \ref{not} below.\\
 Now we want to see that $|A|^{-1}(1)$ is a sufficiently tight obstacle to prevent  $h_{d,\ve}^{n-1}$ from touching $\Sigma$ for $\ve \ra 0$ showing that the auxiliary
distance term in $\pounds_{\ve}$ is dispensable.

\begin{proposition}\label{not}
 \quad For every $d$ there is a neighborhood $V_d$ of $\Sigma$ such that for $\ve \ra 0$ $h_{d,\ve}^{n-1}$ converges to a minimal hypersurface
$h_d^{n-1}$ with obstacle  $|A|^{-1}(d)$  and $h_d^{n-1} \cap V_d = \emptyset$. Moreover, there is a $\beta^1_n > 0$ independent of $H$ such that for some large $d_H$: $ V_d
\supset U_{\beta^1_n / d} (\Sigma) $ the $\beta^1_n / d$-distance tube, for $d \ge d_H$.
\end{proposition}


\textbf{Proof} \quad We start with the case of a minimal cone $C$ singular only in $0$: since $N^{n-1} = \p B_1(0) \cap C$ and $|A|^2$ are analytic we know that
$|A|^{-1}(0) \subset N^{n-1}$ is lower dimensional. \\
We want see that the size of the slices $\p B_\rho(0) \cap |A|^{-1}([0,1])$ of the \emph{funnel set} $|A|^{-1}([0,1])$ shrinks faster than any open subcone for $\rho \ra 0$ and this
prevents an area minimizer with obstacle
$|A|^{-1}(1)$ from entering this funnel too deeply:\\

Since $0$ is the minimum of $a(\omega)^2$ its Taylor expansion in smooth points in $a(\omega)^{-1}(0)$ starts only with second or  higher even order terms. Choosing geodesic
or harmonic coordinates, we may assume that these terms locally uniformly dominate zeroth and first order terms to any desired
extend: since the leading terms are all of order between $2$ and $2k$, for some possibly large but finite $k$.\\

Note that $|A|^{-1}(0) $ is a subcone and along each ray $\lambda^{-1} \cdot |A|(x)=|A|(\lambda \cdot x) $. Hence for $\lambda \in (0,1]$ we have for some $a_{\varsigma} \nearrow
1, b_{\varsigma} \searrow 1$ when $\varsigma \ra 0$ :
\[ a_{\varsigma} \cdot \lambda \cdot \sqrt[2]{\lambda} \cdot  dist(\p
B_1 \cap (|A|^2)^{-1}(\varsigma), \p B_1  \cap (|A|^2)^{-1}(0)) \le \]
\[ dist(\p B_\lambda  \cap (|A|^2)^{-1}(\varsigma), \p B_\lambda  \cap (|A|^2)^{-1}(0)) \le \]
\[ b_{\varsigma} \cdot \lambda \cdot
\sqrt[2k]{\lambda} \cdot  dist(\p B_1 \cap (|A|^2)^{-1}(\varsigma), \p B_1 \cap (|A|^2)^{-1}(0))\]\\
Therefore, when $\lambda \ra 0$, $\ch^{n-1}(\p B_\lambda \cap (|A|^2)^{-1}([0,\varsigma]))$ decreases faster than $\ch^{n-1}(B_\lambda \cap (|A|^2)^{-1}(\varsigma))$. Hence for
small $\ve$ where the obstacle will converge to $|A|^{-1}(0)$ the area minimizer $h_{d,\ve}^{n-1}$  will not exceed a certain shell $\p B_\lambda \cap
(|A|^2)^{-1}([0,\varsigma])$: \\

Since the size of the funnel shrinks faster than linear there is no free area minimizer reaching $0$ and the area of the obstacle $(|A|^2)^{-1}(d)$ is undercut by shells $\p
B_\lambda \cap (|A|^2)^{-1}([0,\varsigma])$. Thus for limit surface $h_d^{n-1}$ with obstacle $(|A|^2)^{-1}(d)$ (which exists since it either coincides with the obstacle or a
subsequence of free minimizers converges in the usual sense (for the transition regions use this for the Plateau problem)) we find an open neighborhood $V_d$ of $0$ with
$h_d^{n-1} \cap V_d = \emptyset$. Obviously this argument reproduces for varying $d$ via scaling invariance of $C$. Therefore we
find a $\beta(C) > 0$ such that $V_d \supset U_{\beta(C) \cdot d} (\Sigma)$.\\
Actually this $\beta$ can be chosen independently of $C$, i.e. $\beta(C)= \beta_n  > 0$: if there is a (converging) sequence of cones $C_m$ in the compact set $\overline {\cal T}_H
\subset C_{n}$ then (since $|A|$ will also converge on smooth parts) we would observe that on compact sets outside $0$ the obstacles converge to the limiting one. But for the
limit cone the minimizer $h_d^{n-1}$ will stay away from $0$. Therefore there is no sequence of cones
$\beta(C_m) \ra 0$.\\

Now we argue as follows: for an analytic minimal hypersurfaces with isolated singularities
this will also be true by cone reduction. \\

In the case of a general singular cone we use induction to get the result: outside $B_1(0)$ small distance tubes and the effective parts of the obstacle around the singular set
$\sigma$ of $C$ become more and more product-like and using the argument above (for $\beta(C_m) \ge const. > 0$) a potential sequence of minimal hypersurfaces
$h_{d,\ve}^{n-1}$  with obstacle containing points $x_{d,\ve}$ with $dist(x_{d,\ve},\sigma) \ra 0$  leads via scalings to
 a local minimizer with obstacle $|A|^{-1}(1)$ (after scaling) on a tangent cone which in this case is a product cone $\R \times C^{n-1}$ and the set $|A|^{-1}(1)$
 is also the product of the set in  $C^{n-1}$ and in $\R$.\\
However the closest minimizer for $C^{n-1}$ taking the product with $\R$ also gives the closest minimizer for $\R \times C^{n-1}$ and thus by definition the obtained minimizer
cannot violate a distance constraint valid for that closest one which in turn is bounded away applying induction.
Again one uses the compactness of the $SC_{n-1}$ to derive that the estimates can be made uniform.\\

Finally we can reduce the result for $H$ to the cone case. Assume for $1/j \ra 0$ $h_{d,1/j}^{n-1} \subset H^n$ contains a sequence of points $y_j$ converging to a point $y_\infty
\in \Sigma$. Then under scaling of $(H,g_H)$ by $(\sqrt{d(y_j,y_\infty)})^{-2}$ the sequence still converges, but at the same time we get an arbitrarily good approximation of the
hypersurface by tangent cones. \qed

 Thus we can directly use $|A|^{-1}(d)$ as an obstacle. Now we want to see how $h_d$ approaches $\Sigma$. On \emph{minimal cones} we note from the flat norm compactness of the space
 of minimal cones that for a minimal cone $C$ any are minimizer  $h_1$ with obstacle $|A|^{-1}(d)$ there are constants $a_1(n), a_2(n) > 0$
 such that $a_1 \le dist(h_1,0)  \le a_2 $. There is a counterpart on arbitrary hypersurfaces

\begin{proposition}\label{pa}
There are constants $k_1(n), k_2(n) > 0$ such that for any $p \in \Sigma$ and $d$ large enough: $k_1 \cdot  1/d \le dist(h_d,p)  \le k_2 \cdot  1/d$.
\end{proposition}

\textbf{Proof} \quad The first claim is covered from the previous Proposition. For the second inequality: $dist(h_d,p)  \le k_2 \cdot  1/d$, we assume that there is no such constant
$k_2
> 0$. For each $d$ we find a $p_d \in \Sigma$ with $dist(h_d,p_d)  \cdot d = sup\{dist(h_d,p)  \cdot d | p \in \Sigma\}$ and for $d \ra \infty$ we get a subsequence converging to some
$p_\infty \in \Sigma$. Following the abstract cone reductions we may assume that $p_d = p_\infty$. Now consider points $q_d \in h_d$ such that $d(q_d,p_d) = dist(h_d,p_d) $.
Next we choose a tangent cone $C$ in $p_\infty$ and notice that for $d \ra \infty$ approximates $d^2 \cdot H$ on arbitrarily large compact regular parts $B_{R_d}(0) \setminus
B_{\varrho_d}(0)) \setminus V_{\xi_d}(\sigma)  \subset C$, i.e.  $R_d \ra \infty, r_d \ra 0, \xi_d \ra 0$ for $d \ra 0$, where $V_{\xi}(\sigma)$ is a cone shaped neighborhood of
$\sigma$ (the cone over the $\xi$-neighborhood of $\p B_1(0) \cap \sigma$). Moreover since $d(q_d, p_d) \ra 0$ we will find (noting the needed scaling by $d^2$) that $d \cdot
d(q_d, p_d) /R_d \ra 0$ and (by assumption)
$d \cdot d(q_d, p_d)  \ra \infty$ (which means that for $d \ra \infty$ increasingly large portions of the $h_d$ are captured in portions of $H$ well-approximated by $C$).\\

Now we use the assumption $d \cdot d(q_d, p_d)  \ra \infty$: we rescale the geometry again by  $d(q_d, p_d)^{-2}$. The effects is that a subsequence of $h_d \cap B_{R_d}(0)
\setminus B_{\varrho_d}(0)) \setminus V_{\xi_d}(\sigma)  \subset C$  converges for $d \ra \infty$ to a non-trivial free area minimizing hypersurface $N^{n-2}$ in $C \setminus
\sigma$ which is complete
within every regular subcone of $C$.\\
From the choice of $p_\infty$ and the relation "$\Sigma \prec \sigma$", between distances to $\Sigma$ versus $\sigma$ we get that $N \subset U_2(\sigma)$ where $U_2$ is
the distance tube of radius $2$ in $C$. However a contraction to $\sigma$ shows there is no such area minimizing $N$.\\

This shows that the $k_2$ exists and does not depend on the base point. However we need the compactness of the space of cones again to conclude that the $k_i$ depend only
on $n$: namely the previous estimate show that $h_d$ scaled by $d^2$ converges (in subsequences) to minimizers with obstacle $|A|^{-1}(1) \subset C$ for any tangent cone. But
here we already noted the existence of uniform estimates which therefore imply that the $k_i$ depend only on $n$.\qed

  Now we define an averaging of $|A|$ and discuss the failure to establish the cone reducibility. However, we will observe that, we still obtain kind of a \emph{weak cone
reducibility}  which suffices for our purposes.\\

Firstly we introduce a \emph{unique} area minimizer ${\cal{J}}$ which will be (called) the \emph{closest} one to $M_1^{m}$: whenever there are two minimizers ${\cal{J}}_1$,
${\cal{J}}_2$ we consider the sets where the supports intersect (and keep this portion) and otherwise for disjoint components we choose the one closer to $M_1^{m}$ and get
another minimizer. Applying this to the family of all minimizers gives a unique
minimizer and henceforth we will mostly consider this minimizer closest  to $M_1^{m}$ and label it by $\hbar$ (with varying suffixes). \\

 Thus, (for $d$ large enough to make this a definition) we will use the neighborhoods of $\Sigma$ made from the set surrounded by $\hbar_d$
 \[W_{d}(\Sigma) := H  \setminus \hbar_d^{n-1} \mbox{
\emph{minus the large volume component of} } H \setminus \hbar_d^{n-1}.\]

 Now we discuss the natural transition properties of these sets $\hbar_d^{n-1}$ between $H$ and its tangent cones.
 The problem becomes visible when we start to define the averaging  $\ap$ of $|A|$ via its level sets as follows: \[x \in \hbar_c \Rightarrow \ap(x) := c.\]
There are two things to worry about: the hypersurfaces $\hbar_k$ for varying $k$ have to be disjoint to ensure that $\ap$ is well-defined (this is obvious from the fact that they
are area minimizers and analytic outside the coincidence sets and the singular set). \\ And secondly,  does every $x \in H \setminus \Sigma$ belong to such a hypersurface ?
Actually, in general, this is not the case: for growing $k$ $\hbar_k$ may \emph{jump}, that is $
W_{l}(\Sigma) \varsubsetneqq  \overline{\bigcap_{k < l} W_{k}(\Sigma)}$. (Since the difference set has positive measure the set of jump levels $\mathfrak{J}_H$ remains countable.)\\

However there is a natural way to extend the definition of  $\ap$ to $H \setminus \Sigma$: for each jump level one  readily checks that (at least) for these levels there are  two
distinct area minimizing hypersurfaces with obstacle $|A|^{-1}(l)$, namely $\hbar_l^{n-1}$ and the limit $h_l^\ast$ of  the $\hbar_k^{n-1}$ for $k \nearrow l$. Since the free parts
(= complement of the coincidence set with the obstacle) of these hypersurfaces are analytic we conclude that there is a least one path component of the free parts of
$\hbar_l^{n-1}$ or $h_l^\ast$ whose interior is entirely disjoint form the other
hypersurface.\\

The space of $H$ between them is not hit by any $\hbar_c^{n-1}$ for a $c \neq l$. Thus for these $x$ we set $\ap(x) := l$. This leads us to the function $\ap$ which is now
well-defined and
continuous everywhere on $H \setminus \Sigma$ and whose level sets do not touch $\Sigma$.\\

 \begin{remark}  We observe that $W_{d}$ can now be seen as a distance tube for the distance measure $1/ \ap$ to $\Sigma$.
 \[W_{d}(\Sigma) = \lev((d,\infty)) = \{x \in H |  1/ \ap < 1/d \}\]
Later on we will often use the complement and thus we give it an own name  \[G_{d}(\Sigma) := H \setminus W_{d}(\Sigma).\]  We  study some of the basic properties of $W_{d}$
(which are also valid for arbitrary area minimizers instead of $\hbar_d^{n-1}$). Obviously $|A| \le d \mbox{ on } G_d$.\\ As an application of the previous arguments applied to the
enhanced singular set cf. sec.12 below we observe that the whole hypersurface $\hbar_d^{n-1}$ shrinks almost proportionally to the distance $1/d$: For each $\ve$ there is a $
\beta^2_n(\ve)
> 0$ independent of $H$ such that almost all of $\hbar_d^{n-1}$ is contained in a distance tube $U_{\beta^2_n / d}$:
\[ (E) \;\;\;\;  {\cal{H}}^{n-1}(\hbar_d^{n-1} \setminus U_{\beta^2_n / d} (\Sigma))\mbox{ / }{\cal{H}}^{n-1}(\hbar_d^{n-1}
\cap U_{\beta^2_n / d} (\Sigma)) \le \ve \; \mbox{ for } \; d \gg 1.  \] Moreover, there is a constant $\kappa =\kappa (M,H,g) > 0$ such that
\[\ch^{n}(W_{d}(\Sigma)) \le \kappa \cdot d^{-8} \;\;  \mbox{ and } \;\; \ch^{n-1}(\hbar_d) \le \kappa \cdot d^{-7}\]
 This is an obvious consequence of Proposition \ref{pa}  (since $\hbar_d$ has a smaller area than the $k_2 \cdot 1/d$-distance tubes modulo the part controlled by (\emph{E})) and
 the fact that ${\cal{H}}^{n-8}(\Sigma) < \infty$  applying the coarea formula.\\
 For use in cone reduction arguments we also notice for products of a minimal hypersurface $N$ (e.g. a cone) with $\R$ and singular set $\Sigma = \Sigma' \times \R$ :
\[W_{d}(\Sigma ) = W_{d}(\Sigma') \times \R \;\; \mbox{ and } \;\; G_{d}(\Sigma) = G_{d}(\Sigma') \times \R\]

 \qed
\end{remark}

Now we resume the discussion concerning the presence of jumps. They clearly cause problems when we pass to tangent cones: on a single cone  the uniqueness of $\hbar_k$
shows that varying $k$ just leads to a rescaling of the $\hbar_k$ and hence \emph{on cones there are no jumps} at all . But there can be different minimizers with the same
obstacle. And if we vary the cone continuously (in flat norm) we also observe that while the level sets of $|A|$ change steadily (on compact smooth parts of the cones) the
hypersurfaces $h_1$ (on
 a family of cones) may again jump in a fashion similar to the jumps described above.\\
In perspective of Proposition \ref{pa} we can infer that  for any $p \in \Sigma$  subsequences $s_i$ of $\hbar_d^{n-1} \subset H$ converge under scaling by $d^2$ to some locally
area
minimizing hypersurfaces $h(s_i)$ with obstacle $|A|^{-1}(1)$ in a tangent cone $C$. However, usually, $h(s_i)$ does \emph{not coincide} with $\hbar_1^{n-1} \subset C$.\\


\emph{Thus we extend the idea and  assign to each level $|A|^{-1}(d) \subset H$ the family $\mathfrak{H}_d$ of all area minimizing hypersurfaces with obstacle $|A|^{-1}(d)$.}\\

Summarizing we observe that the assignment $H \mapsto \mathfrak{H}_d,  d \in \R^{>0}$, is a cone reducible functor. Proposition \ref{pa} shows that the elements in this family
stay
uniformly together in those places which transfer to the tangent cones.\\
Our interest in this paper will be focussed on the growth rate near infinity of certain functions defined on $\lev((0,1)) \subset C$ on tangent cones and this does \emph{not}
depend on the representative of $\mathfrak{H}_1$. Repeated in formal terms $\lev(c) \subset H$, $c \in \R^{>0} \setminus \mathfrak{J}_H$ is as a distinguished representative
of equivalence
classes whose elements can be used compute certain invariants (e.g. the growth rates) and we observe these invariants do \emph{not} depend on the representative. \\
Thus for better readability we argue modulo $\mathfrak{H}_d$ and will henceforth \emph{assume}  $\ap$ is cone reducible and to suppress generic choices we also
\emph{assume} $ \mathfrak{J}_H = \emptyset$.


\vspace{1cm}

\setcounter{section}{4}
\renewcommand{\thesubsection}{\thesection}
\subsection{Positive solutions on $H$}

\bigskip

Now we will construct global conformal deformations $w^{4/n-2} \cdot g$ of $H \setminus \Sigma$ to get metrics with $scal(w^{4/n-2} \cdot g) > 0$ and well controllable
geometry close to $\Sigma$.  To this end we will use Perron solutions. In this section we will prove the existence and uniqueness of such a solution on $H$. We will see that on $H$
(quite different from cones)
the Perron property is actually always satisfied.\\

The existence comes from taking a subsequence of the first eigenfunctions for Dirichlet problems for the singular conformal Laplacian on a sequence of regular domains $K_m
\subset K_{m+1}$ in $H \setminus \Sigma$ with $\bigcup_m K_m = H \setminus \Sigma$. Here we take $ K_m := \lev((0,m)) = G_{m}(\Sigma)$ which will become important when
we want to understand
the limiting behavior of  $\wp_H$ near $\Sigma$.\\

Choosing $\ap$ as a distance measure provides coherence (has the same scaling properties as $|A|$) with the modification of the conformal Laplacian:  We use $ | A | $ as a
weight to get $L_H
= - | A | ^{-2} \cdot ( \triangle + \frac{n-2}{4 (n-1)} \cdot scal_H )$ which is necessary to keep the information concerning the eigenvalue while analyzing via cone reduction.\\

(As before we may assume that $|A|^2$ is \emph{analytic} and the set \emph{$|A|^{-1}(0)$ has $(n-1)$-dimensional measure zero}.)\\

\begin{lemma}\label{exi}
 \quad \[ \lambda_H := \inf_{f \not\equiv 0,smooth, \supp f \subset H \setminus \Sigma}
\frac {\int_{H \setminus \Sigma} |\nabla f|^2 + \frac{n-2}{4 (n-1)} \scal_H f^2} {\int_{H \setminus \Sigma} | A | ^2  \cdot f^2} > 1/10 \] and we can find a smooth function positive
(although not-integrable) function $u_0$ on $H \setminus \Sigma$ with
\[ - \triangle u_0 + \frac{n-2}{4 (n-1)}  \scal_H u_0 = \lambda_H \cdot  | A | ^2  \cdot u_0\]
\end{lemma}

\textbf{Proof} \quad The stability inequality $(A2)$ and $scal_M  > 0$ provide us with the following estimate:
\[ \int_H | \nabla f |^2 + \frac{n-2}{4 (n-1)} scal_H f^2 d A \ge \]
\[\int_H \frac{n}{2 (n-1)} | \nabla f |^2 + \frac{n- 2}{2 (n-1)} f^2 \left( | A |^2 + scal_M
\right) d A  \ge \int_H  \frac{n- 2}{2 (n-1)} | A |^2 f^2 d A \]\\

which gives the estimate for $\lambda_H$. The weight as well as the underlying space are singular and thus we cannot handle $\lambda_H$ as a first eigenvalue
with a corresponding eigenfunction by standard means. \\

But we can construct such a smooth function $u_0$ approximating the problem by a sequence of regular ones. \\

For any $\ve >0$ we find a unique first Dirichlet eigenfunction $u_{m,\ve}$ satisfying
\[ - \triangle u_{m,\ve} + \frac{n-2}{4 (n-1)}\scal_H u_{m,\ve} = \lambda_{m,\ve} \cdot \pounds_\ve^2 \cdot u_{m,\ve}, \quad \lambda_{m,\ve} > 0 \]
with $u_{m,\ve} >0$ on $\op{int}K_m$, $u_{m,\ve} \equiv 0$ on $\p K_m$ and $\int_{B_0} \pounds_\ve^2 \cdot u_{m,\ve}^2 = 1$ for a fixed ball $B_0 \subset H
\setminus \Sigma$.\\

 Since the function space grows for increasing $m$, the eigenvalue $\lambda_{m,\ve}$ decreases monotonically as $m \to \infty$ and hence there is a unique
limit $\lambda_{\infty,\ve} = \lim_{m\to\infty} \lambda_{m,\ve} \ge 0$. Also note that $\lambda_{\infty,\ve} = \lambda_\ve$
\[ \lambda_\ve := \inf_{f \not\equiv 0,smooth, \supp f \subset H \setminus \Sigma} \frac {\int_{H \setminus \Sigma} |\nabla f|^2 + \frac{n-2}{4 (n-1)} \scal_H f^2}
{\int_{H \setminus \Sigma} \pounds_\ve^2  \cdot f^2}  \] since for any function $f$ with compact support in $H \setminus \Sigma$ we eventually have $\supp f
\subset K_m$ for sufficiently large $m$.
\\

{\bf Claim} \quad There is a subsequence of $(u_{m,\ve})_m$ that converges in $C^k$ (for any $k$) to a (not necessarily integrable) limit function $u_\ve >0$ on $H \setminus
\Sigma$ satisfying
\[ - \triangle u_\ve + \frac{n-2}{4 (n-1)} \scal_H \cdot u_\ve = \lambda_\ve \cdot \pounds_\ve^2 \cdot u_\ve. \]
(Note that, unlike $\lambda_\ve$, this limit function {\it may} depend on the
choice of $K_m$.) \\

{\bf Proof} \quad This is a standard application of elliptic estimates and Harnack inequalities. Since such arguments will appear several times later on and the smoothed weight
$\pounds^2_\ve (x)$ might appear unusual, we carry them out
in some detail here. \\

First of all, notice that $\lambda_{m,\ve} \to \lambda_\ve \ge 0$ implies that there exists $c_1 >0$ such that $0 \le \lambda_{m,\ve} \le c_1$ for all $m$. Thus, on every ball $B$ with
compact closure in $H \setminus \Sigma$ the equations
\[ - \triangle u_{m,\ve} + \left( \scal_H - \lambda_{m,\ve} \cdot \pounds_\ve^2 \right) \cdot u_{m,\ve} = 0 \]
have uniformly (in $m$) bounded coefficients. Therefore, we get uniform constants in the interior elliptic estimates
\[ | u_{m,\ve} |_{C^l(B')} \le c_l(B,B') \cdot | u_{m,\ve} |_{L^2(B)} \]
(the $L^2$- and $L_\ve^2$-norms are locally equivalent) and the Harnack inequalities
\[ \sup_{B'} u_{m,\ve} \le \bar c(B,B') \cdot \inf_{B'} u_{m,\ve} \]
for all balls $B' \subset \! \subset B \subset \! \subset H \setminus \Sigma$. \\

Thus, on $B_0$, the $L_\ve^2$-bound = 1 and Harnack's inequality imply upper and lower bounds
\[ c_2'(B_0) > \sup_{B_0} u_{m,\ve} \ge \inf_{B_0} u_{m,\ve} > c_2(B_0) >0 \]
and therefore on a slightly larger ball $B_0' \supset \! \supset B_0$
\[ \sup_{B_0'} u_{m,\ve} \le  c_3 \cdot \inf_{B_0'} u_{m,\ve} \le c_3 \cdot c_2'(B_0), \]
i.e., there is a uniform $L^2$-bound on $B_0'$ and thus a $C^l$-bound on $B_0$ and we may assume that $u_{m,\ve}$ converges in $C^l$ on $B_0$. The limit satisfies $u_\ve \ge
c_0(B_0)
>0$ and the equation
\[ - \triangle u_\ve + \left( \scal_H - \lambda_\ve \cdot \pounds_\ve^2 \right) \cdot u_\ve = 0. \]
\qed

Now is $H \setminus \Sigma$ is connected and hence for any point $x \in H \setminus \Sigma$ outside $B_0$ we can argue by choosing a smooth path $\gamma:[0,1] \to H
\setminus \Sigma$, $\gamma(0) \in B_0$, $\gamma(1)=x$ covered by finitely many overlapping balls $B_1,\dots,B_k$ in order to get $L^2$-estimates: say $B_0 \cap B_1 \neq
\emptyset$; then
\[ \tilde c^{-1} \cdot \sup_{B_1} u_{m,\ve} \le \inf_{B_1} u_{m,\ve} \le \inf_{B_0 \cap B_1} u_{m,\ve} \le \]
\[ \sup_{B_0 \cap B_1} u_{m,\ve} \le \sup_{B_0} u_{m,\ve} \le c_3 \cdot \inf_{B_0} u_{m,\ve} \le c_3 \cdot c_2'(B_0). \]\\
Arguing as for $B_0$ we get a further positively lower and upper bounded  subsequence converging on $B_0 \cup B_1$ and, proceeding by induction, a subsequence
converging in $C^k$ to a limit function $u_\ve >0$ on all of $H \setminus \Sigma$. \\
Next we observe that $\lambda_\ve \ra \lambda_H$ for $\ve \ra 0$ and choosing suitable multiples we may assume that $\int_{B_0} u_{\ve}^2 =1$ for every $\ve > 0$. Thus we can
argue similarily as before and find a $C^k$-converging sequence $ u_{\ve_i}$ for some sequence $\ve_i \ra 0$, $i \ra \infty$ with smooth limit $u_0 > 0$ on $H \setminus \Sigma$
satisfying
\[ - \triangle u_0 + \frac{n-2}{4 (n-1)} \cdot scal_H  \cdot  u_0 = \lambda_H \cdot  | A | ^2  \cdot u_0\] \qed

\begin{remark} At first sight surprisingly one could use decreased eigenvalues: any value $\lambda <  \lambda_H$ could be obtained as an eigenvalue for a positive eigenfunction. (Actually this provides
us
with some degree of freedom to use geometric arguments to derive some estimates later on.)\\

We use two methods to accomplish this decrease: both cases rely on a local $scal$-decreasing deformation as described in [L2]. In the first case  we decrease the scalar
curvature in each step close to $\p K_m$ for the exhausting sequence $K_m$ such that the previously sketched
construction leads to smaller eigenvalues on each $K_m$ and also in the limit on $H \setminus \Sigma$.\\
Thus in the limit we recover the original geometry on $H$ and get as smooth $u_{\lambda} > 0$ with
\[ - \triangle u_{\lambda} + \frac{n-2}{4 (n-1)}  \scal_H u_{\lambda} = \lambda \cdot  | A | ^2  \cdot u_{\lambda}\]
(Clearly, this would have been impossible on a closed manifold.) However the construction has the drawback that the approximating sequence of Dirichlet
solutions lives on domains deformed near $\p K_m$.\\
On the other hand when we allow changes of the final metric on $H$ we can fix one ball $B$ with $\bar{B} \subset H \setminus \Sigma$ and decrease $scal$ on this
ball such that the Dirichlet eigenvalues and that for the limit function decrease to any given extend.\\
The advantage of this approach is that the geometry near $\p K_m$ is now the original one and $\lambda$ can again be characterized by \[ \lambda = \inf_{f \not\equiv 0,smooth,
\supp f \subset H \setminus \Sigma} \frac {\int_{H \setminus \Sigma} |\nabla f|^2 + \frac{n-2}{4 (n-1)} \scal_H f^2} {\int_{H \setminus \Sigma} | A | ^2  \cdot f^2}\] \emph{with
respect to the metric deformed on $B$} (i.e. $\nabla f$, $\scal_H$, the norms and the volume element have to be taken with respect to the new metric whereas $| A |$ remains
unchanged)\qed
\end{remark}

 Now we want to see that (up to multiples) $u_0 $ from \ref{exi} is the unique positive solution of this eigenvalue equation and it has the Perron property. For the latter point it will
 be helpful to see where a Perron solution comes from, that is we want to mimic the classical Perron process. The first point to note is that one cannot apply the standard Perron strategy since the sign of the
linear term is just the converse of the case where one
could apply the maximum principle (cf. [GT], p. 103)). Actually it is a Fredholm alternative argument that helps us to imitate the Perron type approach. \\

To actually begin with, we fix a smoothly bounded neighborhood $V \subset U_\delta(\Sigma) \setminus \Sigma$ of the singular set $\Sigma$ of $H$ within a $\delta$-distance
tube $U_\delta(\Sigma)$. Choosing $\delta \ll 1$ means that $scal_H|_V$ is almost negative: since $ scal_H =  scal_M -
 2 Ric_M (\nu, \nu) - | A |^2 $ the scalar curvature is uniformly upper bounded everywhere and since $\delta \ll 1$ means that in (eventually) most places $|A| \gg 1$
we can scale the whole setting keeping $scal \ll -1$ in most places while $scal_H|_V \ll 1$ everywhere. This can readily be turned into a quantitative statement using tangent cones
where $scal \le 0$ and the zero set is lower dimensional. On $V$ we want to find the smallest solution $u
> 0 $ of the equation (LO)
\[    \triangle u + (\lambda_H |A|^2- \frac{n-2}{4 (n-1)} \cdot scal_H ) \cdot u = 0 \mbox{ with } u \equiv u_{\lambda_H } \mbox{ on } \p V. \]

\begin{lemma}\label{compuniq}
 \quad Let $G \subset H$ be a smoothly bounded domain with $\overline{G} \subset H \setminus \Sigma$. Then the problem
\[ \triangle u + (\lambda_H |A|^2 - \frac{n-2}{4 (n-1)} \cdot scal_H ) \cdot  u = 0 \mbox{ on } \inn{G}  \mbox{ and }  u = \varphi \mbox{ on } \p G \]
has a unique solution for any (continuous) function $\varphi : \p G \to \R$.
\end{lemma}

{\bf Proof} \quad We argue applying the Fredholm alternative ([GT], p.107) for the elliptic operator
\[ \triangle u + (\lambda_H |A|^2 - \frac{n-2}{4 (n-1)} \cdot scal_H ) \cdot  u  \mbox{ on } G \]
That is we show that
\[ \triangle u + (\lambda_H |A|^2 - \frac{n-2}{4 (n-1)} \cdot scal_H ) \cdot  u = 0  \mbox{ on } G \mbox{\:\:\: and \:\:\:}  u \equiv 0 \mbox{ on } \p G\]
has \emph{only} the trivial solution.\\
Otherwise extending a non-trivial solution by $0$ on $H \setminus G$ we get a $w \in H^{1,2}(H \setminus \Sigma)$ with $\mbox{supp } w \subset \bar{G}$ and $\|w\|_{L^2(H
\setminus \Sigma)} = 1$ such that $\int_{H \setminus \Sigma} |\nabla w|^2 + \frac{n-2}{4 (n-1)} \scal_H w^2 / \int_{H \setminus \Sigma} | A | ^2 \cdot w^2 = \lambda_H <
\lambda_H$. But this contradicts the definition of $\lambda_H$.  \qed

This statement is no more true for connected \emph{unbounded} domains. The Martin boundary at infinity can be seen as a measure for the (usually extreme) deviation. We will
need a good
understanding of this boundary to link solutions one $H$ with those on its tangent cones. \\

Here we use Lemma \ref{compuniq} in the study of Perron families on $V$. As usual we call a function $v:V \to \R$ \emph{supersolution} of $\triangle u + (\lambda_H |A|^2- \scal)
u = 0$ if, for any ball $B \subset \! \subset V$ and any solution $u$  on $B$ with $u|_{\p B} \le v|_{\p B}$, it
follows that $u|_B \le v|_B$.\\
In order to ensure that we have got a sufficiently rich class of supersolutions we first notice that the \emph{minimum} of two supersolutions is also a supersolution. Another
operation (a local upgrading of a super- to an actual positive solution) within this class is the {\it lift} $\bar u$ on $B$ of a supersolution $u:V\to\R$, which is defined as follows:
on $B$ we let $\bar u$ be the unique solution of $ \triangle \bar u + (\lambda_H |A|^2- \frac{n-2}{4 (n-1)} \cdot scal_H ) \bar u = 0$ with $\bar u|_{\p B} = u|_{\p B}$ and $\bar u =
u$ on
$V \setminus B$. \\

\begin{lemma} \quad The lift $\bar u$ of a positive supersolution $u$  is still a positive supersolution.
\end{lemma}

{\bf Proof} \quad Choose some ball $B'$ and consider a solution $h$ of equation (LO) with $h \le \bar u$ on $\p B'$. Since $u$ was a supersolution we have $\bar u \le u$ on $B'$ and
thus $h \le u$ on $\p B'$. Hence $h \le u$ on $B'$ and $h \le \bar u$ on $B' \setminus B$ and thus $h \le \bar u$ on $\p (B' \cap B)$. Now if there is point $p \in interior(B' \cap B)$
where $h(p) > \bar u(p)$ then one can take $\max\{h- \bar u, 0\}$ extended by $0$ on $H \setminus (B' \setminus B)$ to get a function $w \in H^{1,2}(H \setminus \Sigma)$ with
$\mbox{supp } w \subset H \setminus \Sigma$ and $\|w\|_{L^2(H \setminus \Sigma)} = 1$ such that $\int_{H \setminus \Sigma} |\nabla w|^2 + \frac{n-2}{4 (n-1)} \scal_H w^2 /
\int_{H \setminus \Sigma} | A | ^2 \cdot w^2 = \lambda^0 < \lambda_H$. But again this contradicts the definition of $\lambda_H$.  Similarly we see that  $\bar u \ge 0$ and by
 Hopf's maximum principle (cf. remark below) $\bar u > 0$. \qed

\begin{remark}
\quad As already mentioned the Hopf's maximum principle applies to general solutions of $\Delta u + g(x) u = 0$ with $g \le 0$, which is precisely \emph{not} our case. However if
u \emph{vanishes} in the point where it is applied one can drop the sign assumption for $g$ (cf. [GT], p.34) and still obtains the critical strict inequality for the \emph{outer
normal derivative} $\p u / \p n > 0$  in an extremal point $q$ of the zero set in the sense that the
interior ball condition for the complement is satisfied (cf. proof of (4.2) for a sample argument) and thus there is a locally (at least relative to this interior ball) unique maximum in $q$.\\
Here and later on we merge this with non-negativity information to utilize this key estimate (from the proof of Hopf's maximum principle) also for our equations.
We just refer to it as the \emph{Hopf's maximum principle}.\\

\end{remark}

After these preliminary considerations, we are now ready to apply the Perron method to our equation. \\To this end, let $S = \{ v:V \to R \,|\, v \mbox{ supersolution}, v > 0, v|_{\p V}
\ge u_{\lambda_H} \}$. $S$ is \emph{non-empty} since at least $u_{\lambda_H} \in S$.

\begin{lemma}
\quad  The function $w(x) = \inf_{v\in S} v(x)$ is positive and satisfies
\[    \triangle w + (\lambda_H |A|^2- \frac{n-2}{4 (n-1)} \cdot scal_H ) \cdot w = 0 \mbox{ on } V \mbox{ with } w \equiv f_H \mbox{ on } \p V. \]
\end{lemma}

{\bf Proof} \quad Obviously $w$ is well defined and non-negative. Let $y$ be an arbitrary point of $V$ and $v_n \in S$ such that $v_n(y) \to w(y)$. By definition, $v_n > 0$ and taking
minima (i.e. replacing $v_n$ by $\min (v_n, v_0)$) we may assume that the sequence $v_n$ is bounded. Now choose a small ball $B \subset V$ around $y$ and consider the lift
$V_n$ of $v_n$ on $B$. We have $V_n \in S$ and therefore $w(y) \le V_n(y) \le v_n (y) \to w(y)$. Moreover, by standard compactness results, we can assume that $V_n$ converges
uniformly on any ball $B' \subset \! \subset B$) to an eigenfunction $v$ on $B$. Clearly $v \ge w$ and $v(y)=w(y)$; we wish to prove that $w=v$ on $B$: So assume there exists $z \in
B$ such that $v(z)
> w(z)$. Choose a function $W \in S$ such that $w(z) \le W(z) < v(z)$ and define $w_k = \min (W, v_k) \in S$ as well as the corresponding lifts
$\Bar w_k$ on $B$. As before we can assume that $\Bar w_k$ converges to an eigenfunction $\bar w$ on $B$ satisfying $w \le \Bar w \le v$ with equality holding at the point $y$.
Hopf's maximum principle gives a contradiction and we conclude that
$v=w$. \\
It remains to show that $w$ is nowhere zero. To see this, choose a point $x_0 \in \p V$  and a small ball $B_R(x_0)$. Let $f_H$ be the unique solution of the equation with
boundary data given by a smooth function $ \phi \ge 0 $ on $(\p B_R(x_0)\cap int V)  \cup  (B_R(x_0) \cap \p  V)$:
\[    \phi \equiv 0 \mbox{ on  } \p B_R(x_0) \cap int V \mbox{ and }  \phi \equiv u \mbox{ near } x_0. \]
Then again by Hopf's maximum principle $f_H > 0$  on $ B_R(x_0)\cap int V$ and since $v \ge f_H$ for every $v\in S$ we have $w > 0$ on $ B_R(x_0)\cap int V$ and joining any
point in $V$ by a chain of balls we analogously get $w > 0$ on $V$. \qed

\begin{corollary}
 \quad The function $w$ (of (3.3)) has the following minimality property: relatively to any neighborhood $W$ of $\Sigma$,
 with $\p W \cap H \setminus \Sigma$ smooth and $W
\subset V$, when $w$ is the smallest positive solution of
 $\triangle \varphi +(\frac{n-2}{4 (n-1)} + \lambda_H) \cdot |A|^2 \cdot \varphi = 0$ on $W$ with $\varphi|_{\p  W } \equiv w$.
\end{corollary}

{\bf Proof} \quad The case where $W = V$ is just a restatement of (3.3). However the minimality still holds for any smaller neighborhood $W$, since otherwise we could take the
minimum of $w$ and a competing smaller solution $w^*$, get a smaller supersolution and via the Perron procedure of (3.3) also a smaller solution with respect to $V$. \qed

Note that the previous arguments work without any compactness assumptions on $H$ or $\Sigma$ and can therefore also be used for tangent cones. We will prove
that $u_{\lambda_H}$ and $w_{C_p}$ coincide with Perron solutions near their singular sets with respect to the boundary data $u_{\lambda_H}$ and $w_{C_p}$.\\

In the compact case we actually have a stronger uniqueness result (as in the smooth closed case). However the remark that follows explains why we stick to Perron solutions.

\begin{lemma}\label{perh}\quad \begin{enumerate}
\item $u_{\lambda_H}$ is a Perron solution.
\item Let $v$ be any solution $> 0$ of $ - \triangle v + \frac{n-2}{4 (n-1)} \cdot scal_H  \cdot  v = \lambda_H \cdot  | A | ^2 \cdot v$ on $H$. Then $v \equiv  u_{\lambda_H}$ up to a
    multiple.
\end{enumerate}
\end{lemma}

 \textbf{Proof} \quad For $(i)$ we note that $u_{\lambda_H}$ is a limit of Dirichlet solutions $u_m$ for the equation on $K_m$
\[ - \triangle u_m + \frac{n-2}{4 (n-1)} \scal_H \cdot u_m = \lambda_{ m , 0} \cdot |A|^2 \cdot u_m. \]
which are also (unique and hence) minimal with respect to the boundary data along $\p W \cup \p K_m$. Since the values along  $\p W$ approach $u_{\lambda_H}$ $C^k$-uniformly
and $u_m \le u_{\lambda_H}$ (since the set of restrictions of positive supersolutions is larger on $K_m$) and from the characterization of Perron solutions
we see that $u_{\lambda_H}$ actually is Perron.\\

For $(ii)$ we take a positive multiple $\mu \cdot v$ such that $\mu \cdot v \ge u_{\lambda_H}$ on the compact set $H \setminus W$. Since $u_{\lambda_H}$ is Perron we infer
that this inequality holds on all of $H$. Now shows the infimum $\mu_0$ of all such $\mu$ and observe again from the Perron property
 that there must be a point $p$ in $H \setminus W$ such that $\mu_0 \cdot v(p) = u_{\lambda_H}(p)$. But Hopf's maximum principle now implies that
$\mu_0 \cdot v \equiv u_{\lambda_H}$.\qed

Now we briefly deduce the fact that $H$ does not carry positive eigenfunctions for eigenvalues $> \lambda_H$\\

\begin{corollary}
 \quad For any $\lambda > \lambda_H$ there does not exist any positive solution of
 $\triangle \varphi +(\frac{n-2}{4 (n-1)} + \lambda) \cdot |A|^2 \cdot \varphi = 0$ on $H$ .\\
\end{corollary}

{\bf Proof} \quad Assume there is a solution $\varphi_\lambda > 0$ of $\triangle \varphi +(\frac{n-2}{4 (n-1)} + \lambda) \cdot |A|^2 \cdot \varphi = 0$ on $H$ for some  $\lambda
> \lambda_H$. The solution $f_H
> 0$ on H corresponding to the eigenvalue $\lambda_H$ is a limit for $m \ra \infty$ and a suitably sequence $\ve_m \ra 0$ of the unique first Dirichlet eigenfunction $u_{m,\ve_m}$
satisfying
\[ - \triangle u_{m,\ve_m} + \frac{n-2}{4 (n-1)}\scal_H u_{m,\ve_m} = \lambda_{m,\ve_m} \cdot \pounds_\ve^2 \cdot u_{m,\ve_m}, \quad \lambda_{m,\ve_m} > 0 \]
with $u_{m,\ve_m} >0$ on $\op{int}K_m$, $u_{m,\ve_m} \equiv 0$ on $\p K_m$ and $\int_{B_0} \pounds_\ve^2 \cdot u_{m,\ve_m}^2 = 1$ for a fixed ball $B_0 \subset H \setminus
\Sigma$ and $\lambda_{m,\ve_m} \ra \lambda_H$. Thus take a large $m \gg 1$, then we find on $K_m$ for a some small $\delta > 0$: $\delta \cdot \varphi_\lambda > u_{m,\ve_m}$
near $\p K_m$ but $\delta \cdot \varphi_\lambda < u_{m,\ve_m}$ on some open set $W$ with regular boundary and $\overline{W} \subset \op{int}K_m$.  But the $\delta \cdot
\varphi_\lambda$ is positive supersolution of $ - \triangle u_{m,\ve_m} + \frac{n-2}{4 (n-1)}\scal_H u_{m,\ve_m} = \lambda_{m,\ve_m} \cdot \pounds_\ve^2 \cdot u_{m,\ve_m},
\quad \lambda_{m,\ve_m} > 0$ on $W$ equal to $u_{m,\ve_m}$ along $\p W$ with $\delta \cdot \varphi_\lambda < u_{m,\ve_m}$. This means that the Perron process gives a
solution $w > 0$, $w \le  \delta \cdot \varphi_\lambda$ on $W$ and thus with $w \neq u_{m,\ve_m}$ but this contradicts (2.2) \qed

\begin{remark} \quad Thus we have a \emph{unique positive first eigenfunction} on $H$ which  indicates that (with some care) we could act as on a smooth compact
 manifold. However a priori this analogy is limited:  beside the fact that we can decrease the eigenvalue the problem on $H \setminus \Sigma$ (where $H$ could also include non-compact cases) has
 many other \emph{non minimal} solutions in the following sense:
 with respect to its boundary data close to $\Sigma$ it is actually \emph{not} true that our
solution would also be a maximizer and thus the uniqueness does not trivialize the minimality statement:\\
Namely, as an instructive sample (which arises as the infinitesimal model in $\Sigma$) we consider an isolated singularity, that is, take a regular cone $C$ and assume for now
that $|A|$ is constant on $\p B_1(0) \cap C$ and that we know that  the equation $- \triangle w + \frac{n-2}{4 (n-1)} \cdot scal_H \cdot  w = \lambda_H \cdot  | A | ^2  \cdot w$
degenerates to an ODE of second order. In this case we observe that for the boundary value $1$ on $\p B_1(0) \cap C$ one gets two canonical linearly independent positive
solutions $r^{-\alpha_i}, \alpha_i > 0, i=1,2$ with different pole order in $0$. The tamer one, let us say $r^{-\alpha_1}$, is a limit of the Dirichlet solutions of regular domains
exhausting $C \setminus 0$ and these solutions can be written as $(1 + \ve) \cdot r^{-\alpha_1}
-\ve \cdot r^{-\alpha_2}$.\\

Note also that this already provides a counterexample to uniqueness even of positive solutions in the non-compact case. Thus it is more adequate to consider the minimality
(towards the singular set) as the distinctive feature.
\end{remark}\qed

\textbf{Induced Solutions on Cones}: Now we consider any tangent cone $C_p$ in a point $p \in \Sigma$ (we may set $p=0 \in \R^n$) and show how $u_{\lambda_H}$ induces a
solution of $ \triangle w + (\frac{n-2}{4 (n-1)} + \lambda_H \cdot  | A | ^2 ) \cdot w = 0$ on $C_p$: Consider $H \cap G_m \setminus G_{m/2} \cap B_{R/m}(p)$ and scale the
metric by $m^2$, then we observe for any $R > 0$ and $a \ge 2$
\[ d_{B_{R}(0)} (m \cdot (H \cap G_{a \cdot m} \setminus G_{m/a}), C_p \cap G_a \setminus G_{1/a}) \to 0 \]
and since $C_p \cap G_a \setminus G_{1/a} \cap B_R(0)$ is contained in the smooth part of  $C_p$ we can infer $C^k$-convergence.\\
With respect to this identification we can consider restrictions
\[w(a,R) := w|_{H \cap G_{a \cdot m}\setminus G_{m/a} \cap B_{R/m}(p)}\]  as being defined on $C_p \cap G_a \setminus G_{1/a} \cap B_R(0)$\\

Since $\Delta w, scal$ and $|A|^2$ scale in the same way $w$ also solves the equation on this domain. Now choose a fixed ball $B
\subset C_p \cap G_a \setminus G_{1/a} \cap B_{R_0}(0)$ for some $R_0 > 0$ and normalize for any $R \ge  R_0$, $w(a,R)$ to $L^2$-norm $1$ on $B$.\\
For $a \ra \infty$ and $R \ra \infty$ (e.g. setting $R = const. \cdot a$ cf. the properties of $G_m$ discussed in the previous section)), one can find a
$C^k$-converging subsequence of $w(a,R)$ with a limit solution $w_{C_p} > 0$ defined on $C_p \setminus \sigma$.\\
Thus after $L^2$-normalizations the solution $u_{\lambda_H}$ defined on $H$ $C^k$-approximates $w_{C_p}$ compactly .\\

In the next few sections we want to show that
 $w_{C_p}$ is a Perron solution.\\

\vspace{1cm}

\setcounter{section}{5}
\renewcommand{\thesubsection}{\thesection}
\subsection{Dirichlet Problems and Positive Solutions on Cones}

\bigskip

We start with some understanding of the space of positive solutions of Dirichlet problems on cones for exhausting families of domains. This provides a good insight into the
space of solutions on cones.

\begin{proposition} \quad There is precisely one positive solution $w_a$ (up to multiples) for the following problem on $G_a \subset C$:
 \[\triangle \varphi +\left(\frac{n-2}{4 (n-1)} + \lambda_H \right)\cdot |A|^2 \cdot \varphi = 0, \; \varphi|_{\lev(\{a\}) } \equiv 0.\]
\end{proposition}

In this section we first analyze two more symmetric cases: the case of product cones and of regular subcones.  Then, in the next section, we also get and use some coarse
information for the behaviour of solutions on $C$. Finally we
compose these details to get a proof of this proposition.\\

As a technical tool we need  boundary Harnack inequalities (also called Carleson inequalities). Those versions existing in the literature (as in [CFMS] or [CS]) apply directly only to
operators (typically $div(A(x)  \cdot \nabla u) = 0$) \emph{without} a zeroth order term (i.e. a term containing only the function but non of its derivatives). Thus we deduce the
following extension for our operator

\begin{lemma} \label{car} \quad  For any solution of the boundary problem \[ \triangle u + (\lambda_H  + \frac{n-2}{4 (n-1)} )  \cdot |A|^2 \cdot  u  = 0 \mbox{ on } G_a \setminus G_1 \subset C
\times \R \] with $u \equiv 0 \mbox{ on } \p G_a \cup \p G_1$ and any point $p \in \p G_a \cup \p G_1$ we have for some $\rho(n,a) \in (0,1/4 \cdot min\{1/2,dist(\p G_a, \p G_1)\})$
for $z \in B_{\rho}(p) \cap G_a \setminus G_1, \rho \in (0,\rho(n,a))$ and some $C(n) > 0$:
\[ u(z) \le C(n) \cdot sup\{u(x) |  x \in B_{\rho^2}(x_p) \cap G_a \setminus G_1\}  \] for some $ x_p \in G_a \setminus G_1 $ with
 $ dist(x_p, \p G_a) = \rho $ resp. $ dist(x_p, \p G_1)  = \rho
$ on $\p B_{\rho}(p) \cap G_a \setminus G_1.$\\
\end{lemma}

{\bf Proof} \quad  From [CFMS] or [CS],Th.11.5  we know that for any nonnegative harmonic function $u \ge 0$, $\Delta u = 0$ on $B_1(0) \cap \R^{n-1} \times R^{\ge 0} $ with $u
\equiv 0$ along $\R^{n-1} \times \{0\}$ we have
\[(\ast) \qquad u(z) \le \kappa(n) \cdot u(1/2 \cdot e_n) \mbox{ for any }   z \in B_{1/2}(0) \cap \R^{n-1} \times R^{\ge 0}   \]
where $e_n$ is the $n-th$ coordinate unit vector and $\kappa(n) > 2$.\\
Now for $p \in \p G_a \cup \p G_1$ (for easier notations we assume  $p \in \p G_1$) and we consider
 \[ \triangle u + (\lambda_H  + \frac{n-2}{4 (n-1)} )  \cdot |A|^2 \cdot  u  = 0 \mbox{ on }
G_a \setminus G_1 \subset C \times \R \]
and consider the pair  $ \p G_1 \subset  G_a \setminus G_1 $  as the substitute for $\R^{n-1} \times \{0\} \subset B_1(0) \cap \R^{n-1} \times R^{\ge 0} $. We will show that we
can actually take $C(n) := 2 \cdot \kappa(n)$: assume there is a sequence $u_m$ of nonnegative solutions defined at least on $B_{2/m}(p) \cap G_a \setminus G_1$ such that
there are points $z_m$ in $B_{1/m}(p) \cap G_a \setminus G_1$ with
\[ u(z_m) \ge  2 \cdot \kappa(n) \cdot sup\{u(x) |  x \in B_{1/m^2}(x_p) \cap G_a \setminus G_1\}  \]
After scaling by $m^2$ we can normalize their $L^2$-norm on $B_{2/m}(p) \cap G_a \setminus G_1$ (new radius = 2) to $1$ and we can argue that the domains converge in
$C^k$-norm to $B_{2}(0) \cap \R^{n-1} \times R^{\ge 0}$ while $|A|^2 \ra 0$ $C^k$-uniformly under this scaling. Thus boundary regularity applied to the family of operators
$\triangle u + \tau \cdot (\lambda_H  + \frac{n-2}{4 (n-1)} )  \cdot |A|^2 \cdot  u  = 0$ for $\tau \in [0,1]$ gives the same upper $C^k$-bound on the scaled ball $B_{2/m}(p) \cap
G_a \setminus G_1$ for all $u_m$ and from the smooth convergence of the domains we infer the existence of a $C^{k-1}$-converging subsequence of $u_m$ on these scaled balls
$B_{2/m}(p) \cap G_a \setminus G_1$ to a nonnegative harmonic function $u_\infty$ on $B_{2}(0) \cap \R^{n-1} \times R^{\ge 0}$ with $L^2$-norm $1$ and $u_\infty|_{\bar
B_{2}(0) \cap \R^{n-1} \times \{ 0\}} \equiv 0$ and a limit point $z_\infty \in  \bar B_1(0) \cap  \R^{n-1} \times R^{\ge 0}$ outside $\R^{n-1} \times \{ 0\} $
\[ u(z_\infty) \ge  2 \cdot \kappa(n) \cdot u_\infty(1/2 \cdot e_n) \] which violates $(\ast)$. \qed

From this and with the aid of Perron solutions  we can understand the Martin boundary for our operator on \emph{difference sets of the} $G_a$'s where we are additionally
interested in the effect of varying the size of this difference in order to get global results on $G_a$. These results could also be seen as an extension of some
Phragm\'en-Lindel\"{o}f type results
in the classical theory of harmonic functions which have been treated by rather different methods cf. [P](8.6) for probabilistic approaches for harmonic functions.\\

We will present a detailed argument for the following result as a sample of how to \emph{utilize} Perron solutions as a transparent alternative for quite a bit more involved
techniques.

\begin{proposition}\label{exp}  \quad  The boundary problem \[ \triangle u + (\lambda_H  + \frac{n-2}{4 (n-1)} )  \cdot |A|^2 \cdot  u  = 0 \mbox{ on } G_a \setminus G_1 \subset C \times \R
\] with $u \equiv 0 \mbox{ on } \p G_a \cup \p G_1$ has two generating solutions (positive on $\inn{G_a \setminus G_1}$)  $\Psi^\pm(x,t) =  h(\pm t) \cdot \psi(x)$ with $h(t) = exp(\gamma_a
\cdot t)$ with $\gamma_a > \gamma_{a'} > 0$
for $1 <a < a'$: every solution $v >0$ on the interior can be written as \[v(x,t) = \alpha_+ \cdot \Psi^+(x,t) + \alpha_- \cdot \Psi^-(x,t)\] for some $\alpha_+,\alpha_- \ge 0$. \\
\end{proposition}

{\bf Proof} \quad Consider the restrictions of the non-trivial solutions to
$ G_a \setminus G_1 \times \{0\}$ whose $L^2$-norm on this set is $1$, we check that the subset of \emph{nonnegative} solutions $\cal{S}^+$ is \emph{compact}:\\
 We get uniform $C^k$ estimates via interior Harnack inequality
 (and elliptic regularity) for $C^k$ norms on  $ G_{a-\ve} \setminus G_{1 + \ve} \times [-1, 1]$ for some given $\ve \in (0,\frac{a-1}{100})$ and let $\varpi > 0$ be
 the constant in the Harnack inequality on $ G_{a-\ve} \setminus G_{1 + \ve} \times [-1, 1] \subset  G_{a-\ve/2} \setminus G_{1 + \ve/1} \times [-2, 2]$:\\
 for any nonnegative solution $v$:   $\sup v \le \varpi \cdot \inf v$ on $ G_{a-\ve} \setminus G_{1 + \ve} \times [-1, 1]$ and since
 $\sup_{G_{a-\ve} \setminus G_{1 + \ve} \times \{0\}} v \le \sup_{G_{a-\ve} \setminus G_{1 + \ve} \times [-1, 1]} v$ and
 $\inf_{G_{a-\ve} \setminus G_{1 + \ve} \times [-1, 1]} v \le \inf_{G_{a-\ve} \setminus G_{1 + \ve} \times \{0\}} v $ the inequality persists under restrictions. Thus
{\footnotesize { \[\varpi^{-2} \cdot (\sup_{G_{a-\ve} \setminus G_{1 + \ve} \times \{0\}} u)^2 \cdot \int_{G_{a-\ve} \setminus G_{1 + \ve} \times \{0\}} dA = \int_{G_{a-\ve}
\setminus G_{1 + \ve} \times \{0\}} \varpi^{-2} \cdot (\sup_{G_{a-\ve} \setminus G_{1 + \ve} \times \{0\}} u) ^2 dA  \]
\[ \le \int_{G_{a-\ve} \setminus G_{1 + \ve} \times \{0\}} (\inf_{G_{a-\ve} \setminus G_{1 + \ve} \times \{0\}} u)^2 dA
\le \int_{G_{a-\ve} \setminus G_{1 + \ve} \times \{0\}} u^2 dA \]}}
Therefore applying the Harnack inequality again we have {\footnotesize { \[\int_{G_{a-\ve} \setminus G_{1 + \ve} \times [-1, 1]} u^2 dV \le \int_{G_{a-\ve} \setminus G_{1 + \ve}
\times [-1, 1]} (\sup_{G_{a-\ve} \setminus G_{1 + \ve} \times [-1, 1]} u)^2  dV \le\]
\[ \varpi^{2} \cdot \int_{G_{a-\ve} \setminus G_{1 + \ve} \times [-1, 1]}
(\inf_{G_{a-\ve} \setminus G_{1 + \ve} \times [-1, 1]} u)^2  dV\le \varpi^{2} \cdot \int_{G_{a-\ve} \setminus G_{1 + \ve} \times [-1, 1]}
(\inf_{G_{a-\ve} \setminus G_{1 + \ve} \times \{0\}} u)^2  dV\]
\[\le \frac{\varpi^{4} \cdot \int_{G_{a-\ve} \setminus G_{1 + \ve} \times [-1, 1]}(
\int_{G_{a-\ve} \setminus G_{1 + \ve} \times \{0\}} u^2 dA)dV}{\int_{G_{a-\ve} \setminus G_{1 + \ve} \times \{0\}} dA } \le \varpi^{4} \cdot \frac{\int_{G_{a-\ve} \setminus G_{1 + \ve}
\times [-1, 1]}dV}{ \int_{G_{a-\ve} \setminus G_{1 + \ve} \times \{0\}} dA} \]}} Thus we get $L^2$-and hence  $C^0$-bounds on $G_{a-\ve} \setminus G_{1 + \ve} \times [-1, 1]$. Using
the Carleson inequality \ref{car}  we can extend these to get uniform $C^0$-estimates and thus $L^2$-estimates on $ G_a \setminus G_1 \times [-1, 1]$. Hence from boundary
regularity we get uniform estimates for $C^k$ norms on $ G_a \setminus G_1 \times [-1, 1]$. Iteratively we get uniform $C^k$-estimates on each  $ G_a \setminus G_1 \times [-m,
m] \setminus (-(m-1), m-1)$, $m \ge 2, m \in \Z$.
In particular any sequence in $\cal{S}^+$ has a compactly $C^{k-1}$-converging subsequence on $ G_a \setminus G_1 \times \R$.\\

The space $\cal{S}^+$ is \emph{non-empty}: consider the (uniquely solvable) boundary value problem on $G_a \setminus G_1 \times[-j, m]$ with boundary value $= 0$ on $\p (G_a
\setminus G_1) \times [-j, m] \cup \inn{G_a} \setminus G_1 \times \{m\}$ and some positive function $\zeta_{j,m}$ on $\inn{G_a} \setminus G_1 \times \{-j\}$. First notice that
the solution $\psi_{j,m}$ is \emph{positive} in the interior:  consider the family of equations
\[ \triangle u + t (\lambda_H  + \frac{n-2}{4 (n-1)} )  \cdot |A|^2 \cdot  u = 0, \]
with the same boundary conditions for $t \in [0,1]$. These problems also have unique solutions $u_t$ and thus $u_t$ depends continuously on $t$. But $f_H > 0$ in the interior by
the minimum principle for harmonic functions and assuming there are $t \in [0,1]$ with $u_t(x) < 0$ for some interior point $x$ we also find a $\tau \in [0,t]$ with $u_\tau
(x_\tau) \ge 0$ but $u_\tau (x_\tau) = 0$ for some interior point $x_\tau$ (since otherwise the normal derivative along the boundary does not vanish anywhere by Hopf's maximum
principle)
but in turn for $u_\tau (x_\tau) = 0$ Hopf's maximum principle leads to a contradiction. Hence (for $t = 1$) we have $\psi_{j,m}$ is \emph{positive} in the interior. \\

Now one can choose constants $\lambda_{j,m} > 0$ such that $\| \lambda_{j,m} \cdot \psi_{j,m}\|_{L^2( G_a \setminus G_1 \times [0, 1])} = 1$ and arguing as above we get,
sending $j \ra \infty$, a compactly $C^k$-converging subsequence with a limit $\psi_{\infty,m}$, \emph{positive} on $ \inn{G_a} \setminus G_1 \times \R^{< m}$, for the boundary
problem
\[ \triangle \psi_{\infty,m} + (\lambda_H  + \frac{n-2}{4 (n-1)} )  \cdot |A|^2 \cdot  \psi_{\infty,m}  = 0 \mbox{ on } G_a \setminus
G_1 \cap C \times \R^{\le m}
\] with $\psi_{\infty,m} \equiv 0 \mbox{ on } \p G_a \cup \p G_1 \cup \inn{G_a} \setminus G_1 \times \{m\}$.
Now we repeat this argument for $\psi_{\infty,m}$: after normalizing the $L^2( G_a \setminus G_1 \times [0, 1])$-norm to $1$ we send  $m \ra \infty$ and get a limit solution
$\psi_{\infty,\infty}$ \emph{positive} on
$ \inn{G_a} \setminus G_1 \times \R$.\\

Knowing that the space of positive solutions is non-empty we can now apply the Perron process to get Perron solutions $\pi_l$, $l \ge, l \in \Z$ for the boundary value problems on
$G_a \setminus G_1 \times \R^{\ge -l}$ with boundary value $= 0$ on $\p (G_a \setminus G_1) \times \R^{\ge -l}$ and $\psi_{\infty,\infty}$ on $\inn{G_a} \setminus G_1 \times
\{-l\}$. Since $\psi_{\infty,\infty} \ge \pi_l \ge 0$ we get for $l \ra \infty$ a compactly $C^k$-converging subsequence with a limit solution $\pi_\infty \ge 0$. Considering the
definitions of these functions and of the Perron solution we observe from the unique solvability of the compactly bordered problems that $\pi_\infty$ is the Perron solution of the
boundary value problems on $G_a \setminus G_1 \times \R^{\ge -l}$ with boundary value $= 0$ on $\p (G_a \setminus G_1) \times \R^{\ge -l}$ and $\pi_\infty$ on $\inn{G_a}
\setminus G_1 \times \{-l\}$
and actually that $\pi_\infty \equiv \psi_{\infty,\infty}$.\\

We claim that the two solutions $\pi^+(x) := \pi_\infty(x)$ and $\pi^-(x) := \pi_\infty(-x)$ (which are linear independent since there is no bounded solution as is seen as from an
obvious localization argument via lemma(3.1) ) generate the space of all solutions: To this end we check that there is a $\lambda^{\pm} > 0$ such that for any $v \in \cal{S}$ we
have: $\lambda^+ \cdot \pi_\infty \le v$ on $ G_a \setminus G_1 \times \R^{\ge 0}$ and that
$\lambda^- \cdot \pi_\infty \ge v$  on $ G_a \setminus G_1 \times \R^{\le 0}$.\\

Hopf's maximum principle shows that the normal derivative of any non-trivial and non-negative solution does not vanish in any point in $\p G_a \cup \p G_1$. Since translations of
solutions in $\cal{S}^+$ along $\R$ lead again to solutions in $\cal{S}^+$ we notice from the compactness of  $\cal{S}^+$ that we have uniform positive upper and lower  bounds
$a > b > 0$ for the normal derivatives in $\p G_a \cup \p G_1 \times \{t\})$ of solutions for any $t \in \R$ when their ${L^2( G_a \setminus G_1 \times \{t\})}$-norm is normalized
to $1$.

Thus we can choose a $\lambda^+ > 0$ such that $\lambda^+ \cdot \pi_\infty \le v$ on $ G_a \setminus G_1 \times \{0\}$ for any $v \in \cal{S}^+$ and the Perron property gives
this inequality on $ G_a \setminus G_1 \times \R^{\ge 0}$.\\  Now assume that there is \emph{no} $\lambda^-$ such that $\lambda^- \cdot \pi_\infty \ge v$
 on $ G_a \setminus G_1 \times \R^{\le 0}$ for any $v \in \cal{S}^+$.
Then we can find a diverging series of slices (i.e.  a sequence of points $t_m \ra -\infty$ while $m \ra \infty$) and a $w \in \cal{S}^+$ with
\[ (1) \;\;\;  m \cdot \pi_\infty(x_m) \le w(x_m) \mbox{ for some } x_m \in \inn{G_a \setminus G_1} \times \{t_m\}.\]
Since the compactness of $\cal{S}^+$ gives us uniform estimates for the normal derivatives at the boundary there is some constant $a > 1$:
\[(2) \;\;\; \frac{1}{a} \cdot \frac{\pi_\infty|_{ G_a \setminus G_1 \times \{t_m\}}}{ |\pi_\infty|_{L^2( G_a \setminus G_1 \times \{t_m\})}}
 \le \frac{w|_{ G_a \setminus G_1 \times \{t_m\}}}{|w|_{L^2( G_a \setminus G_1 \times \{t_m\})}}  \le
 a \cdot \frac{\pi_\infty|_{ G_a \setminus G_1 \times \{t_m\}}}{ |\pi_\infty|_{L^2( G_a \setminus G_1 \times \{t_m\})}}\]

Combining $(1)$ and $(2)$ there are suitably large $0 \ll k \ll l$ and a constant $c > 0$ such that on the interior of the slices:
\[c \cdot  \pi_\infty|_{\inn{G_a \setminus G_1} \times \{t_k\}} >
w|_{\inn{G_a \setminus G_1} \times \{t_k\}} \;\mbox{\; and \;}\; c \cdot  \pi_\infty|_{\inn{G_a \setminus G_1} \times \{t_l\}}  <
w|_{\inn{G_a \setminus G_1} \times \{t_l\}}\]
But this contradicts the Perron property of $\pi_\infty$ on $ G_a \setminus G_1 \times \R^{\ge t_l}$.\\

Now we choose any  $f \in \cal{S}^+$ and form  $F(x) :=  f(x) + f(-x)$. We claim $F \equiv c' \cdot (\pi^+(x) + \pi^-(x))$ for some $c' > 0$: Consider the \emph{infimum} $\lambda_H$
of all  $\lambda > 0$ with  $F \le \lambda \cdot (\pi^+(x) + \pi^-(x))$. The previous discussion shows that the set of these $\lambda$ is nonempty. If $F(x) < \lambda_H \cdot
(\pi^+(x) + \pi^-(x))$ in some interior point Hopf's maximum principle says the same on the whole interior and  the outward normal derivatives of $ \lambda_H \cdot (\pi^+(x) +
\pi^-(x)) - F(x)$ are negative everywhere. Thus for some tiny $\ve > 0$ we
still have $F(x) < (\lambda_H - \ve) \cdot (\pi^+(x) + \pi^-(x))$ on $ G_a \setminus G_1 \times [-j,j]$ for some given $j$.\\
But for $t \ra \infty$ (applying again Hopf's maximum principle) we observe on $G_a \setminus G_1 \times \{t\}$: $F/ \lambda_H \cdot (\pi^+(x) + \pi^-(x)) \ra 1$ and thus on the
interior $(\lambda_H - \ve) \cdot (\pi^+(x) + \pi^-(x)) < F(x)$. But repeating the argument
for the existence of  $\lambda^-$ such that $\lambda^- \cdot \pi_\infty \ge v$ this also contradicts the Perron property of $\pi_\infty$.\\

In particular $f \equiv c' \cdot \pi^+(x)$ on the  slice $G_a \setminus G_1 \times \{0\}$. But if we translate $\pi^+(x) + \pi^-(x)$ by $t$ and use the same argument for any slice
$G_a \setminus G_1 \times \{t\}$ (not just for $t=0$) we find that $f(t) \equiv c_t' \cdot \pi^+(x)$ for a suitable $c_t'> 0$. In other words every nonnegative solution of the
equation is of the form $h(t) \cdot \pi^+(x)$ and inserting gives $(3)$:
\[\frac{\p^2 h(t)}{\p t^2} \cdot  \pi^+(x) + h(t) \cdot \triangle_{G_a \setminus G_1 \times \{0\}} \pi^+(x) +
(\lambda_H  + \frac{n-2}{4 (n-1)} )  \cdot |A|^2 \cdot  h(t) \cdot \pi^+(x)  = 0 \]
Since $h(t)$ does not depend on $x$ we evaluate this in the interior (where $\pi^+(x) > 0$): this is a linear second order ODE with constant coefficients $\frac{\p^2 h(t)}{\p t^2} +
C \cdot h(t) = 0$. Since we already know that $h(t) > 0$ and unbounded we conclude $C < 0$ and get the solutions
\[h(t) = \alpha_1 \cdot exp(\sqrt{(-C)} \cdot t) + \alpha_2 \cdot exp(-\sqrt{(-C)} \cdot t), \;\;\;\; \alpha_1,\alpha_2 > 0\]

Finally we want to see how $C$ depends on $a$:  we claim  $\sqrt{(-C(a))} > \sqrt{(-C({a'}))} > 0$ for $1 <a < a'$. To this end
 consider $\alpha_1 =\alpha_2 = 1/2$ (thus $h(0) = 1$) and note $\frac{\p^2 h(t)}{\p t^2} = -C =|C|> 0$:
\[  -\triangle_{G_a \setminus G_1 \times \{0\}} \pi^+(x) -
(\lambda_H  + \frac{n-2}{4 (n-1)} )  \cdot |A|^2  \cdot \pi^+(x)  = |C| \cdot \pi^+(x)\] Thus the solution can be described as  the \emph{first} eigenfunction (since it is positive) of
the Dirichlet problem for the operator $\triangle_{G_a \setminus G_1 \times \{0\}} + (\lambda_H  + \frac{n-2}{4 (n-1)} )  \cdot |A|^2$ being the minimizer of the variational
integral
\[ \int_{G_a \setminus G_1 \times \{0\}} |\nabla f|^2 - (\lambda_H  + \frac{n-2}{4 (n-1)} ) \cdot |A|^2 \cdot  f^2\] with over all
$f \in H_0^{1,2}(G_a \setminus G_1 \times \{0\}), |f|_{L^2}=1$.  We infer via Hopf's maximum principle that the solution is uniquely determined (up to a multiple) and we observe
that the infimum $|C|(a)$ of the variational integral strictly decreases for increasing $a$:\\ This is easily seen by constructing a test function for $a' > a$: take a tiny tubular
neighborhood $U$ of $\p G_a$ in $G_{a'} \setminus G_1$ such that $U$ can be identified with $\p G_a \times (-\delta,\delta)$ for some $\delta \ll 1$ (oriented such that $\p G_a
\times \{-\delta\} \subset G_a \setminus G_1 $). On this set substitute $\pi_a^+(z,s)$, $(z,s) \in \p G_a \times (-\delta,\delta) \subset G_a \setminus G_1$ extended by zero on
$\p G_a \times (0,\delta)$ for $\pi_a^+(z,1/2 \cdot s - 1/2 \cdot \delta)$.\\ Thus in the variational integral the \emph{integral} over the $|\nabla f|^2$-term will decrease linearly
whereas that over the $(\lambda_H  + \frac{n-2}{4 (n-1)} ) \cdot |A|^2 \cdot  f^2$-term grows linearly in $\delta$. \qed

Now we turn to \emph{regular subcones} of a  given cone $C$ with singular set $\sigma \subset C$: this means we choose a domain $G \subset \p B_1(0) \cap C$ with compact
closure and smooth boundary  such that $\overline{G} \subset \p B_1(0) \cap C \setminus \sigma$ and consider the  \[\mbox{\emph{cone over G}: \quad}  C(G) := \{x = t \cdot z \,
|  \,  t \ge 0, z \in \overline{G}\} \subset C\]

\begin{proposition}\label{con} \quad  The boundary problem \[ \triangle u + (\lambda_H  + \frac{n-2}{4 (n-1)} )  \cdot |A|^2 \cdot  u  = 0 \mbox{ on } \inn{C(G)} \subset C
\] with $u \equiv 0 \mbox{ on } \p C(G) \setminus \{0\}$ has two generating solutions (positive on $\inn{C(G)} $)  $\Psi^\pm_G(x,t) =  c_G(\omega) \cdot r^{\alpha^\pm_G}$,
 $\alpha_\pm = - \frac{n-2}{2} \pm \sqrt{ \left( \frac{n-2}{2} \right)^2 - \mu_G}$, where $\mu_G$ is the first eigenvalue for the Dirichlet problem on $G$.
\end{proposition}

{\bf Proof} \quad Due to the scaling invariant setting and the transformation properties of the entities within the equation we have for such a solution $u$ that for any $b > 0$
$u_b := u(b \cdot x)$ solves \[ \triangle u + (\lambda_H  + \frac{n-2}{4 (n-1)} )  \cdot |A|^2 \cdot  u = 0  \mbox{ on } G_{a \cdot b} \setminus G_b
\]
From that on the proof is completely similar to the previous one and leads to analogous separation of variables. Next we derive the claimed form for the exponent.\\

 $c(\omega)$ solves {\footnotesize \[ (CW)    \;\;\;
    \left(\alpha^2 + (n-2) \alpha \right) \cdot c(\omega) + \left( \triangle_{C \cap \p B_1(0)}  + \left( \frac{n-2}{4 (n-1)}+\lambda_H\right) a(\omega)^2 \right) c(\omega) = 0\]} and is a
    Perron solution.

 $c(\omega)$ is an eigenfunction of the operator $\triangle_S+( \frac{n-2}{4 (n-1)}+\lambda_H) a(\omega)^2$ with eigenvalue
$\mu=-\alpha^2-(n-2)\alpha$, and therefore
\[  \alpha_\pm = - \frac{n-2}{2} \pm \sqrt{ \left( \frac{n-2}{2} \right)^2 - \mu}. \]
Since  $c(\omega) > 0$ $\mu$ is the first eigenvalue for the Dirichlet problem on $G$. \qed

\begin{corollary}\label{corsingcone} \quad  Let $D_i \subset \p B_1(0) \cap C$ be a sequence of domains with $D_i \subset D_{i+1}$, $\bigcup_i D_i = \p B_1(0) \cap C$.\\ Then we have for
sufficiently large $i$:  $\mu_{D_i} > 0$ and thus $\alpha _{D_i}\in (- \frac{n-2}{2}, 0)$.
\end{corollary}

{\bf Proof} \quad To understand the Dirichlet eigenvalue equation on ${D_i}$ \[\left(\alpha^2 + (n-2) \alpha \right) \cdot c(\omega) + \left( \triangle_{C \cap \p B_1(0)}  + \left(
\frac{n-2}{4 (n-1)}+\lambda_H\right) a(\omega)^2 \right) c(\omega) = 0\] we consider the variational integral \[ \inf \{ \int_{D_i} | \nabla f |^2 - \left( \frac{n-2}{4
(n-1)}+\lambda_H\right) a(\omega)^2 \cdot  f^2 d A \;|\; f \in H_0^{1,2}({D_i}), |f|_{L^2} = 1 \}\] and want to show that it becomes negative when $i$ is large enough. Since $- \left(
\frac{n-2}{4 (n-1)}+\lambda_H\right) a(\omega)^2 <0$ almost everywhere, we want to find a function which $1$ in the interior of ${D_i}$  except for a small tube ${U_i}$ around
$\p {D_i} $ where it falls off to $0$. However
this part produces positive contributions from $| \nabla f |^2$. Thus our goal is to define such functions  $f_i$ with $\int_{U_i} | \nabla f_i |^2 \ra 0$ for $i \ra \infty$.\\

Although this is not an inductive argument it may be helpful to see the case of isolated singularities first: assume a ball surrounding a singular point is well approximated by a
ball in a regular tangent cone $C$: here we take $U_i = B_{2^{-i}} \setminus B_{2^{-(i+1)}}$ and $f_i(x) = F_i(r)$ where $F_i \in C^{\infty}(\R,[0,1])$ with $F_i = 0$ on
$\R^{<2^{-(i+1)}}$, $F_i = 1$ on $\R^{>2^{-i}}$. Since the cone has dim $C > 2$ we see that $F_i$ can be chosen such that $\int_{U_i} |
\nabla f_i |^2 \le c_n \cdot Vol(B_{2^{-i}} \setminus B_{2^{-(i+1)}}) \cdot (2^{-i})^{-2} \le c'_n \cdot (2^{-i})^{7}  / (2^{-i})^2\ra 0$ for $i \ra \infty$.\\

Now for higher dimensional singularities we start with a finite covering by balls $B_{r_i}(p_i)$ around singular points $p_i \in \Sigma$ which are well-approximated by balls in
appropriated tangent cones $C_i$ and chosen according to the definition of the $n-3$-dimensional Hausdorff measure. Note the $n-7$-dimensional measure of $\Sigma$ is
already finite and hence all
higher dimensional Hausdorff measures vanish. But we actually only need that the codimension is $> 2$: we can assume $\sum_j r_j^{n-2} < \ve$ for any prescribed $\ve > 0$.\\
Then we construct the desired cut-off function via induction on each of the balls $B_{r_j}(0) \subset C_i$ separately  and take the product of the function extended by $1$ to the
remainder of $D_i$ and by $0$ close to the singular set.\\ Then  using $\int_{C \cap \p B_1(0)}| \nabla \prod_j f_j |^2 \le \sum_j \int_{C_i} |\nabla f_j |^2 \le  \sum_j r_j ^{n} \cdot r_j
^{-2}  < \ve$ we get for large $i$ for $D_i$
\[ \alpha^2 + (n-2) \alpha  < 0 \; \mbox{ and hence } \; \mu >0 \]
which implies $\alpha \in (- \frac{n-2}{2}, 0)$.\qed

Let $C$ be a regular tangent cone in some point $p \in \Sigma$. Writing points in $C$ in polar coordinates, i.e. $(\omega,r) \in C$ where $r \ge $ is the distance to the tip and
$\omega$ a point in $\p B_1(0) \cap C$ we have

\begin{corollary}
\quad The space of positive solutions on $C$ is spanned by  two linear independent positive solutions $f^\pm_C$ on $C$.\\ Every positive solution admits a separation of variables,
more precisely, we get $f^\pm_C = c(\omega) \cdot r^{\alpha_\pm}$ for some positive function $c(\omega) > 0$, $\alpha_- < \alpha_+ <0$ solving {\small
\[(CW) \;\; \left(\alpha_\pm^2 + (n-2) \alpha_\pm \right) \cdot c(\omega) + \left( \triangle_{\p B_1(0) \cap C} + \left( \frac{n-2}{4 (n-1)}+\lambda_H\right) a(\omega)^2 \right) c(\omega)
= 0 \] } The Perron solution satisfies $\wp_C = f^+_C = c(\omega) \cdot r^{\alpha_+}$with
$\alpha_+ \in (- \frac{n-2}{2},0)$.\\
Moreover for $n=7$-dimensional singular minimal cones there are uniform bounds $- \frac{n-2}{2} < \theta_1 < \theta_2 < 0$
such that for any  cone $\alpha_+ \in ( \theta_1,\theta_2 )$.\\
\end{corollary}

{\bf Proof} \quad The fact that there are precisely these two solutions resp. that the exponents satisfy these estimates follows  directly from the arguments in Proposition
\ref{con} resp. Corollary \ref{corsingcone} above. The Perron property is part of the argument that shows that there are just these two generating functions. \qed

Finally we turn to the complementary case $C \setminus C(G)$, that is a cone shaped neighborhood of $\sigma$. Here the space of positive solutions is \emph{infinite
dimensional}: the point is that the regular part of  $C \setminus C(G) \cap \p B_1(0) $ is now an open set and we may apply the techniques of sec.3 to decrease the eigenvalue
$\mu_ {C \setminus C(G)}$ in \ref{con} which is compensated from a differently chosen $\alpha_G$ which allows us to find a large family of positive solutions.

\begin{corollary}\label{open} \quad  The boundary problem \[ \triangle u + (\lambda_H  + \frac{n-2}{4 (n-1)} )  \cdot |A|^2 \cdot  u  = 0 \mbox{ on } \inn{C \setminus C(G)} \subset C
\] with $u \equiv 0 \mbox{ on } \p C(G) \setminus \{0\}$ has an infinite set of linear independent solutions (positive on $\inn{C \setminus C(G)} \subset C$)
$\Psi_\mu(x,r) = \psi_\mu(x) \cdot r^{\alpha_\mu}$, $\alpha_\mu \in (-\infty, \alpha_G)$.\\
Moreover when $D_i \subset \p B_1(0) \cap C$ is a sequence of domains with $D_i \subset \overline{D_i} \subset D_{i+1}$, $\bigcup_i D_i = \p B_1(0) \cap C$: \[\alpha_{\p B_1(0)
\cap C \setminus D_i} \ra \infty \mbox{\;\; for \;\;} i \ra \infty.\].
\end{corollary}

{\bf Proof} \quad For $i \ra \infty$ the domain $\p B_1(0) \cap C \setminus D_i$ shrinks to $\p B_1(0) \cap \sigma$. Thus the eigenvalue diverges.\qed

A simple but clearly important side effect is

\begin{corollary}
\quad On a \textbf{singular} minimal cone $C$ there are infinitely many linear independent positive solutions.
\end{corollary}

Thus our aim will be to characterize a particular solution which will be the Perron solution.

\vspace{1cm}

\setcounter{section}{6}
\renewcommand{\thesubsection}{\thesection}
\subsection{Uniqueness results}

Using the previous special cases of strips and subcones we can now prove the uniqueness of positive Dirichlet solutions on the (non-compact) regular subdomains $G_a \subset
C$. In terms of Martin theory this means that the Martin boundary (at infinity) is a single point. \\
This will be used to derive the heredity principle leading from $\wp_H$ to $\wp_C$.\\

The strategy is to show that assuming there are two positive solutions which are linear independent we can concentrate their deviation in a small almost cylindrical tube around
$\p G_a \subset C$ and then we get a contradiction from the potential shape of the solution on these tubes when we compare them on tangent cones.

\begin{proposition}
 \quad There is  a unique positive solution $w_a$ (up to multiples) for the following problem on $G_a \subset C$:
 \[\triangle \varphi +(\frac{n-2}{4 (n-1)} + \lambda_H) \cdot |A|^2 \cdot \varphi = 0, \; \varphi|_{\p  G_a} \equiv 0.\]
\end{proposition}

 We start with a little but useful observation: positive solutions of $\triangle \varphi +(\frac{n-2}{4 (n-1)} + \lambda_H) \cdot |A|^2 \cdot \varphi = 0$ on cones clearly satisfy a Harnack
inequality on any given pair of balls $B \subset \! \subset B'$. However what is special about this situation is that the constant
 in the Harnack inequality ($\sup u \le c \cdot \inf u$) is \emph{scaling invariant}: with any solution $v(x) > 0$ on $B'$ we also have a solution $v(1/\tau \cdot x) > 0$ on $\tau \cdot B'$, $\tau > 0$ with the
\emph{same} (optimal) constant in the respective Harnack inequalities as is seen from the direct transition between solutions on $B'$ and on  the scaled copy $\tau \cdot B'$.\\

To get an idea of the shape of such solutions on products we start with a proof of a weaker version of the Proposition for product cones
(which actually asserts $w(x,t) = W(x)$, i.e. $w$ is translation invariant):\\

\begin{lemma}
 \quad For any  positive solution $w$ for the following problem on $G_a \subset C \times \R$:
 \[\triangle \varphi +(\frac{n-2}{4 (n-1)} + \lambda_H) \cdot |A|^2 \cdot \varphi = 0, \; \varphi|_{\p  G_a} \equiv 0.\]
we find that $w_a(x,t) / exp(t) \ra 0$ for $t \ra \pm\infty$ where $(x,t) \in C \times \R$.
\end{lemma}

{\bf Proof} \quad  We first prove the following\\

\textbf{Claim}\quad  Otherwise we observe
\[h_1(t) \cdot \psi_1(x) \le w_a(x,t) \le h_2(t) \cdot \psi_2(x) \mbox{ with }\]
\[h_i(t) = \alpha_1 \cdot exp(\sqrt{\kappa_i} \cdot t) + \alpha_2 \cdot exp(-\sqrt{\kappa_i} \cdot t),
\mbox { with } \kappa_2 >  \kappa_1 > 0 \mbox { and } \alpha_1,\alpha_2 \ge 0\]
such that $\kappa_1,\alpha_1,\alpha_2$ depend on $w_a$ while $\kappa_2$ is the same for any positive solution. \\

{\bf Proof} \quad The upper bound comes from the following general consideration: First consider tubes $T_i \subset G_2$ around $ \p G_1 \times \{i\}$, $i \in \Z^{\ge 0}$  with
$T_{i+1} = J(T_i)$ where $J :C \times \R \ra C \times \R$ is the translation by $1$: $J(x,t)=(x,t+1)$ such that $T_i \cap T_{i+1} \neq \emptyset$ and note that for any positive
solution $\phi$ via Harnack inequalities (which hold for some fixed $C \ge 1$):
\[\inf_{T_0}\phi \ge C \cdot \sup_{T_0}\phi \ge C \cdot \inf_{T_1}\phi \ge C^2 \cdot \sup_{T_1}\phi \ge C^2 \cdot \inf_{T_2}\phi \ge C^3 \cdot \sup_{T_2}\phi \ge \dots\]
Thus we have, let us say, on the whole strip $G_{0.9} \setminus G_{1.1}$ a growth estimate: $\phi(x,t) \le exp(k \cdot t) \cdot \inf_{\p G_1 \times \{0\}}\phi$ for some $k \ge 1$. But
this also shows that the growth along any strip  $G_{\beta} \setminus G_{\beta+1}$ is upper bounded by $c_\beta \cdot exp(k \cdot t)$ for the \emph{same} $k$: namely, we
extend the growth estimates along $G_{0.9} \setminus G_{1.1}$ appending a chain of balls $B_r(x,0),..,B_r(y,0)$ from say $(x,0) \in G_{\beta} \setminus G_{\beta+1}$ to $(y,0) \in
G_{0.9} \setminus G_{1.1}$, then we follow $ G_{0.9} \setminus G_{1.1}$ until we reach $(y,T) \in G_{0.9} \setminus G_{1.1}$ for any  $T \gg 0$, exit the track $ G_{0.9} \setminus
G_{1.1}$ (at this station $T$) and return to $(x,T)$ using the same chain of balls translated by $T$. Then we get (for this $x \in C$) a fixed constant $c_x > 0$ such that $\phi(x,T)
\le c_x \cdot \phi(x,0) \cdot exp(k \cdot T)$ and this holds for any $T \in \R$ since we always append the same two fixed (up to translation) finite sets of balls to the chain of balls
along $G_{0.9} \setminus G_{1.1}$. Since $C \times \R$ becomes locally nearly flat near infinity we could actually get approximatively the same $k$ as in the Euclidean space and
we realize that $\kappa_2$ is independent of $w_a$,
 $a$ and even of $C$.\qed

Now for $a' \gg a$ notice that $w_{a'}(x,t) :=  w(\frac{a'}{a} \cdot x, \frac{a'}{a} \cdot t)$ solves the equation on $G_{a'}$. However in this latter case the exponent of $exp$ will
grow from $\kappa_1$ to $\frac{a'}{a} \cdot \kappa_1$ while $\kappa_2$ remains unchanged and
 this eventually violates the upper bound.\\

\qed

Now let us assume we have two linearly independent positive solutions $\phi$, $\Upsilon$ on $G_a \subset C$ with vanishing boundary value:.\\

\begin{lemma}\label{harnack}
 Along any  regular subcone  $C(G)$ we have a constant $\kappa > 0$ such that for $r \ge 1$ and $(r,\omega) \in \R^{> 0} \times
G$: \[ \phi(r,\omega) \le  \kappa \cdot \Upsilon(r,\omega) \:\: \mbox{ or } \:\: \Upsilon(r,\omega) \le \kappa \cdot \phi(r,\omega)\].
\end{lemma}

{\bf Proof} \quad  Otherwise we have a ray $\R^{> 0} \times \{\omega \}$ with $\omega \in G$ and some increasing sequence $r_i \ra \infty$ such that
\[\phi(r_{2i},\omega)/\Upsilon(r_{2i},\omega) \ra \infty \:\: \mbox{ whereas } \:\: \phi(r_{2i+1},\omega)/\Upsilon(r_{2i+1},\omega) \ra 0\].
With some function $v(x)$ we know that $v(\lambda_H\cdot x)$, $\lambda_H> 0$ also solves $\triangle \varphi +(\frac{n-2}{4 (n-1)} + \lambda_H) \cdot |A|^2 \cdot \varphi = 0$.
Thus we can also say that we have two sequences of positive solutions $\phi_i(x) := \phi(r_{2i} \cdot x)$ and  $\Upsilon_i(x) := \Upsilon(r_{2i} \cdot x)$ such that
\[\phi_i(1,\omega)/\Upsilon_i(1,\omega) \ra \infty \:\: \mbox{ whereas } \]
{\small{\[\phi_i(r_{2i-1}/r_{2i},\omega)/\Upsilon_i(r_{2i-1}/r_{2i},\omega) \ra 0 \mbox{ and } \phi_i(r_{2i+1}/r_{2i},\omega)/\Upsilon_i(r_{2i+1}/r_{2i},\omega) \ra 0\]}}

Using the scaling invariant Harnack inequality we infer that $r_{2i-1}/r_{2i} \ra 0$ and $r_{2i+1}/r_{2i} \ra \infty$ and for any domain $U$ with compact closure in $\overline{U}
\subset C \setminus \sigma$ with $(1,\omega) \in U$ we find uniformly \[\phi_i(x)/\Upsilon_i(x) \ra \infty \:\: \mbox{ for } x \in U \]

But that implies that for some exhausting sequence of such domains $U_i \subset U_{i+1}$ with $\bigcup_i U_i = C \setminus \sigma$ $\Psi_i(x) := \phi_i(x) - \Upsilon_i(x) > 0$ on
$U_i$ and such that on the endpoints $x_i, y_i$ of $\p U_i \cap \R^{> 0} \times \{\omega \}$  $\Psi_i(x_i) = \Psi_i(y_i) = 0$ (of course we may assume that $U_i \cap \R^{> 0}
\times \{\omega \}$ is connected). Now let $z_i \in U_i \cap \R^{> 0} \times \{\omega \}$ be chosen such that $\Psi_i(z_i) = max_{x \in U_i \cap \R^{> 0} \times \{\omega
\}}\Psi_i(x)$ and consider the functions $\Phi_i(r,\varpi) := \Psi(\|z_i\| \cdot r, \varpi)/\Psi(\|z_i\|, \omega)$ then we get a subsequence of compactly $C^k$-converging positive
functions on domains $V_i$ that still (via Harnack inequality) exhaust $C \setminus \sigma$ and get a positive limit solution $\Phi_\infty(r,\varpi)$ on $C \setminus \sigma$ such
$\Phi_\infty(1,\omega) = 1$ and
$\Phi_\infty(r,\omega) \le 1$. The scaling invariance of the Harnack inequality shows that  $\Phi_\infty(r,\varpi)$ remains bounded along each ray   $\R^{> 0} \times \{\varpi \}$.\\

 Now recall
the shape of positive solutions $u_G$ on a regular subcone  $C(G) \subset C$: $c(\omega)  \cdot (a \cdot r^{\alpha^-} + b \cdot  r^{\alpha^-})$ for some $a,b \ge 0$ and
$\alpha^-<\alpha^+< 0$. In particular we note again that on $C(G) \cap \overline{B_1}(0)$ the Perron solution $w$ for the boundary data $c(\omega)$ along $C(G) \cap \p B_1(0)$
and (although this is a void condition) $0$ along $\p C(G) \cap B_1(0)$ has a pole $r^{\alpha^+}$. But for some suitable $c > 0$ we note that $c \cdot \Phi_\infty(r,\varpi) > w$
along $\p (C(G) \cap \overline{B_1}(0))$ although $c \cdot \Phi_\infty(r,\varpi) $ is bounded contradicting the fact that $w$ is unbounded. \qed

Now we want to see that such inequalities still hold on $G_a \subset C.$\\

We start with asymptotic growth estimates for any given positive solution $\phi$ of the described problem: solutions on regular subcones $C(G)$ decay slower than $\phi$ near
infinity:
\begin{lemma}\label{up}
 \quad Let $\chi$ be a positive solution for the following problem on $C(G) \cap G_a \subset C$:
 \[\triangle \varphi +(\frac{n-2}{4 (n-1)} + \lambda_H) \cdot |A|^2 \cdot \varphi = 0, \; \varphi|_{\p  (C(G) \cap G_a)} \equiv 0.\]
Then for any $c > 0$ with $\phi < c \cdot \chi$ in some $p \in C(G) \cap G_a$  there is an unbounded domain $D \subset  C(G) \cap G_a$, $p \in D$ with
\[\phi < c \cdot \chi \mbox{ on } D\]
\end{lemma}

{\bf Proof} \quad Per definition $\phi > 0$ while $\chi = 0$ along $\p  (C(G) \cap G_a) \setminus G_a$. Thus if the statement does not hold we can find (via Sard's Lemma) a
constant $k
>0$ such that $\phi < k \cdot \chi $ on a \emph{bounded} domain $U \subset C(G) \cap G_a$ with smooth boundary and $\phi =k \cdot \chi $ along
$\p U$. But there is no positive solution with vanishing boundary data on $U$ acc. Lemma \ref{compuniq}.\qed

Completely similar we can derive asymptotic growth estimates by comparison with solutions on $C  \setminus C(G)$ :
\begin{lemma}\label{low}
 \quad Let $\chi$ be a positive solution on $C  \setminus C(G)$ of
 \[\triangle \varphi +(\frac{n-2}{4 (n-1)} + \lambda_H) \cdot |A|^2 \cdot \varphi = 0,  \; \varphi|_{\p  C(G) } \equiv 0\]
Then for any $c > 0$ with $\chi < c \cdot \phi$ in some $p \in C \setminus C(G) \cap G_a$  there is an unbounded domain $D \subset  C \setminus  C(G) \cap G_a$, $p \in D$ with
\[\chi < c \cdot \phi \mbox{ on } D\]
\end{lemma}

The \emph{proof} is almost literally the same as that of the previous lemma.

\begin{lemma}
For some finite $\kappa_0 > 0$ we have
\[ \phi(r,\omega) \le  \kappa_0 \cdot \Upsilon(r,\omega) \:\: \mbox{ or } \:\: \Upsilon(r,\omega) \le \kappa_0 \cdot \phi(r,\omega)\]
on $G_a \subset C.$
\end{lemma}

{\bf Proof} \quad  According to the previous lemma we may assume, that on any regular subcone $ \phi(r,\omega) \le  \kappa \cdot \Upsilon(r,\omega)$ for a suitable $\kappa$.
Thus when we assume that such an inequality does not hold on $G_a \subset C$ we can find a sequence of regular cones $C_m \subset C$, $C_m \subset C_{m+1}$, $\bigcup_m
C_m = C$ and
a sequence of points $p_m \in G_a \setminus C_m$ with $ \phi(p_m) > m \cdot \Upsilon(p_m)$ and $\phi(r,\omega) - m \cdot \Upsilon(r,\omega) < 0$ on $G_a  \cap C_m$. \\

Thus we conclude from Lemma \ref{compuniq} that there is an unbounded domain $D_m \subset G_a) \setminus C_m$ with $p_m \in D_m$ such that  $ \phi(p_m) > m \cdot
\Upsilon(p_m)$ on $D_m$ and   $ \phi(p_m) = m \cdot \Upsilon(p_m)$ on $\p D_m$.\\

The  growth of $f_m:= \phi(p_m) - m \cdot \Upsilon(p_m)$ on $D_m$ exceeds a given polynomial order $\ge r^k$ for suitably large $m$ as is seen from Lemma \ref{up} and
Corollary \ref{corsingcone}.  We may assume that $f_m(p_m) =1$ and maximal on the slice of points in $D_m \cap \p B_{|p_m|}(0)$ and consider the sequence of pointed spaces
$(D_m, p_m)$ and consider the limit. Note that this is \emph{not} a scaling argument: running along branches of $G_a$ to infinity means that (subsequences of) these pointed
spaces converge to $G_a $ within some tangent cone which is a product cone. \\ The choice of $p_m$ leads to two possibilities: if $p_m$ stays within a upper bounded distance to
$\p G_a$. Then we observe exponential growth in directions parallel to $\p G_a$ acc. Proposition \ref{exp} but since it defines a function on $G_a $ on the tangent cone this
cannot not exist: just consider two rays with different distance to $\p G_a$ then the exponential growth rate must coincide as is seen from a Harnack inequality (which works
with the same constants along these rays). On the other hand Proposition \ref{exp} shows that the growth rate depends on the size of the strip which could be chosen arbitrarily
in this course of this argument. Thus this case will not show up.\\ Thus we may assume that $dist(p_m, G_a) \ra \infty$, but then we get via rescaling (by $\sqrt{ dist(p_m,
G_a)}$) positive solutions on the tangent cone which is bounded in radial direction (which does not exist as was shown in Lemma \ref{harnack}). \qed

Noting that this applies to any two positive solutions we can use the previous result to actually derive

\begin{lemma}\label{evr} For some finite $\kappa_1 > 0$ we have
\[ \phi(r,\omega) \equiv \kappa_1 \cdot \Upsilon(r,\omega) \]
on $G_a \subset C.$
\end{lemma}

{\bf Proof} \quad  We choose $ \kappa_1 := \inf \{ \kappa \,|\, \phi(r,\omega) \le  \kappa \cdot \Upsilon(r,\omega)   \}$. Assume that $\phi(r,\omega) \neq \kappa_1 \cdot
\Upsilon(r,\omega)$, then Hopf's maximum principle shows that actually $\phi(r,\omega) < \kappa_1 \cdot \Upsilon(r,\omega)$ on $G_a \subset C$. However now we can consider
the positive solution $\psi(r,\omega) := \kappa_1 \cdot \Upsilon(r,\omega) - \phi(r,\omega)$. Thus, for some  finite $\kappa_2 > 0$ we have
\[ \phi(r,\omega) \le  \kappa_2 \cdot \psi(r,\omega) \:\: \mbox{ or } \:\: \psi(r,\omega) \le \kappa_2 \cdot \phi(r,\omega)\]
on $G_a \subset C$. The definition of $\kappa_1$ implies that $\psi(r,\omega) \le \kappa_2 \cdot \phi(r,\omega)$\\
But this means for some $\kappa_3 > 0$: $ \Upsilon(r,\omega) \le \kappa_3 \cdot \phi(r,\omega)$. Henceforth we may assume $\kappa_3$ is chosen minimal and that this
inequality is \emph{strict}
 (again via Hopf's maximum principle). Thus we first note that we can sharpen the previous lemma to: for some finite $\kappa_4
> 0$ we have
\[ \phi(r,\omega) \le  \kappa_4 \cdot \Upsilon(r,\omega) \:\: \mbox{ \emph{and }} \:\: \Upsilon(r,\omega) \le \kappa^{-1}_4 \cdot \phi(r,\omega)\]
on $G_a \subset C$. But since such a relation is now true for any two positive solutions we also get a $\kappa_5 > 0$ with
\[ \kappa_5 \cdot \phi(r,\omega) \le   \kappa_3 \cdot \phi(r,\omega)-\Upsilon(r,\omega) \le \kappa^{-1}_5 \cdot \phi(r,\omega).\]
However this is a contradiction to the minimality of $\kappa_3 $ and hence
\[ \phi(r,\omega) \equiv \kappa_1 \cdot \Upsilon(r,\omega) .\] \qed

\vspace{1cm}

\setcounter{section}{7}
\renewcommand{\thesubsection}{\thesection}
\subsection{Perron solutions on tangent cones}

Now we will study Perron solutions on an arbitrary \emph{singular} cone  $C$. There are several ways to obtain (and thus characterize)  $\wp_C$ also showing that  $\wp_C$ is
amenable to a separation of variables.

\begin{proposition}\label{per} There is an (up to multiples) unique  Perron solution $\wp_C$.
\begin{enumerate}
\item $\wp_C$  is the $C^k$-compact limit of a (suitably normalized) subsequence of Perron solutions $\wp_{C(G_\rho)}  = c_\rho(\omega) \cdot r^{\alpha_\rho}$ on
    $C(G_\rho)$, where $G_\rho$ is the complement of the $\rho$-neighborhood of $\sigma \cap \p B_1(0)$  for $\rho \ra 0$.
\item $\wp_C  = c(\omega) \cdot r^{\alpha}$ where $\alpha = \lim_{\rho \ra 0} \alpha_\rho$. \item $c(\omega)$ solves {\footnotesize \[ (CW)    \;\;\; \left(\alpha^2 + (n-2) \alpha
    \right) \cdot c(\omega) + \left( \triangle_{C \cap \p B_1(0)}  + \left( \frac{n-2}{4 (n-1)}+\lambda_H\right) a(\omega)^2 \right) c(\omega) = 0\]} and is a Perron solution.
    \item $\wp_C$  is the $C^k$-compact limit of a (suitably normalized) subsequence of $w_a$ on $G_a$.
    \item There are constants $-\frac{n-2}{2}  < \theta_1(n) < \theta_2(n) < 0$ with $\alpha \in (\theta_1(n),\theta_2(n))$.
\end{enumerate}
\end{proposition}

{\bf Proof}  \quad $(i) - (iii)$: The existence proof showing that the $C^k$-compact limit of a (suitably normalized) subsequence of Perron solutions on $C({G_\rho})$
$\Psi^+_{G_\rho}(x,t) = c_{G_\rho}(\omega) \cdot r^{\alpha^+_{G_\rho}}$,
 $\alpha_+ = - \frac{n-2}{2} + \sqrt{ \left( \frac{n-2}{2} \right)^2 - \mu_{G_\rho}}$, where $\mu_{G_\rho}$ is the first eigenvalue for the Dirichlet problem on ${G_\rho}$ for
$\rho \ra 0$ is actually a Perron solution  $\wp_C  = c(\omega) \cdot r^{\alpha}$ uses that  $c(\omega)$ is Perron: the equation for $\wp_C$ becomes
\[ \frac{\p^2 \wp_C}{\p r^2} + \frac{n-1}{r} \frac{\p \wp_C}{\p r} + \frac{1}{r^2} \triangle_{C \cap \p B_1(0)} \wp_C + \left( \frac{n-2}{4 (n-1)} + \lambda_H\right) \frac{a(\omega)^2}{r^2} \wp_C = 0, \]
where $\triangle_{C \cap \p B_1(0)} $ denotes the Laplacian on $C \cap \p B_1(0)$ and thus, substituting $\wp_C = c(\omega) \cdot r^\alpha$, we see that the equation has a
form that allows to repeat the arguments from Lemma \ref{perh}. (By induction $ c(\omega)$ also converges to $\infty$ in $C \cap \p B_1(0) \cap \sigma$ which allows us to
ignore $ \left(\alpha^2 + (n-2) \alpha \right) \cdot c(\omega)$.)
\[  \left(\alpha^2 + (n-2) \alpha \right) \cdot c(\omega) + \left(\triangle_{C \cap \p B_1(0)}  + \left( \frac{n-2}{4 (n-1)}+\lambda_H\right) a(\omega)^2 \right) c(\omega) = 0. \]
Now the Perron property of  $c(\omega)$ can be seen completely similar as in  Lemma  \ref{perh}. From that we will consider neighborhoods $V_\ve$ of $\p B_1(0) \cap \sigma$ in
$\p B_1(0) \cap C$ and take full dimensional neighborhoods $C(V_\ve) \cap B_{R+\ve}(0) \setminus B_{R-\ve}(0)$ to readily check that $\wp_C $ has the Perron property for all
points in $\sigma \setminus \{0\}$ with respect to this neighborhood.
(For $0$ we can just take a distance ball.)\\
$(v)$ follows from Corollary \ref{corsingcone} and the compactness of the space of singular cones.\\

Next, we prove that $\wp_C$ is the  \emph{unique }Perron solution (up to multiples). Let $v > 0$ be another Perron solution, then (the argument of) \ref{up} and the scaling
invariant Harnack inequality show that along any ray in $C$ $v$ has the same asymptotic growth rate as $\wp_C$. In particular near infinity it converges to zero. As in
\ref{harnack} we infer that, say,  $v \le c \cdot \wp_C$, for some $c_G
> 0$ for any regular subcone on $C(G) \cap (C \setminus B_1(0))$. If this inequality does not hold everywhere on $C(G) \cap (C \setminus B_1(0))$ we argue  acc. \ref{low} that we get a
contradiction to the  existence of solutions acc. \ref{open} that converge to infinity. But then we can turn around (now considering $B_r(0) \cap C$  and argue via Perron property
that the
inequality also holds everywhere on $C$. As in \ref{evr} we conclude that $\wp_C$ and $v$ are equal up to a multiple.\\

Finally we check $(iv)$: To see that $\wp_C$  is the $C^k$-compact limit of a (suitably normalized) subsequence of $w_a$ on $G_a$ it suffices according to the previous argument
to see that this limit has the same asymptotic growth rate near
infinity as $\wp_C$.\\

We choose some fixed ball $B \subset \overline{B} \subset G_1 \subset C \setminus \sigma$ normalize  $w_a$ to $|w_a|_{L^2(B)}=1$, for $a \ge 1$  and compare them with the
Perron solutions ($L^2$-normalized on $B$) on $C(D_i)$, where $D_i \subset \p B_1(0) \cap C$ is a sequence of domains with $D_i \subset \overline{D_i} \subset D_{i+1}$,
$\bigcup_i D_i = \p B_1(0) \cap C$ for $i \ra \infty$. Running through the $i$'s we see that \ref{up} implies that for each $a$  $w_a$ has the same asymptotic growth rate near
infinity as $\wp_C$.  \qed

The estimate  $\alpha  < 0$ means that the cone will become in a sense acuter when being deformed with
 $c(\omega) \cdot r^\alpha$ moreover there are some important geometric properties:

\begin{lemma} \quad   Let $C \subset \R^n$ be a cone and $g$ the induced metric on $C$. Then $C$ equipped with the metric $\tilde g := (c(\omega) r^\alpha)^{4/n-2} \cdot g$ is again a
\textbf{cone} (although not embed) with finite distance between $0$ and any other point of $C$:\\ $(C,\tilde g)$ is isometric to any of copy scaled around $0$ and can be
reparametrized as $c(\omega)^{4/n-2}  \cdot g_{\R} + r^2 \cdot g_{\p B_1(0) \cap C}$ and the scalar curvature in a point (on the ray $(\omega,t) \in C$) with new distance $\rho$
to $0$ is equal
to $\frac{4 (n-1)}{2 |\alpha| } \cdot \lambda_H \cdot c(\omega)^{4 (n-3)/n-2} \cdot a(\omega)^2/\rho^2$. \\
\end{lemma}

{\bf Proof} \quad We first note that $(r^\alpha)^{2/n-2}$ for $- \frac{n-2}{2} < \alpha$ is integrable on $\R^+$ and therefore the distance between $0$ and any other point of $C$
remains finite. Moreover, a conformal deformation with a function of the form $r^\beta$ gives a conical metric on $C$ again: the radius resp. the distance of $\p B_t(0)$ to $0$ in
the transformed metric is $t^{1 + \alpha \cdot \frac{2}{n-2}}$   resp.  $ \frac{n-2}{2 |\alpha| } \cdot  t^{1 + \alpha \cdot \frac{2}{n-2}}$.

Thus, we can restrict to deformations only depending on $\omega$: $\tilde g = c'(\omega) \cdot g$.

Let $f: C \to C$ be  the map $(r,\omega) \mapsto (\lambda_Hr, \omega)$. Of course $g$ has the property $ f^*g = \lambda^2 \cdot g$
but in fact, $\tilde g$ shares this property:
\begin{eqnarray*}
(f^* \tilde g)_{(r,\omega)} (X,Y) &=& c' (f(r,\omega)) \cdot g_{(f(r,\omega))} (df_{(r,\omega)} X, df_{(r,\omega)} Y) \\
&=& c'(\omega) \cdot f^*g (X,Y) = \lambda^2 \; \tilde g(X,Y).
\end{eqnarray*}
Also the cone shape of $C$ gives $|A|(\omega,r) = a(\omega)/r$
\[ \lambda_H \cdot  a(\omega)^2/r^2   =  \frac{n-2}{4 (n-1)} scal_{(c(\omega) r^\alpha)^{4/n-2} \cdot g_H}(\omega,r) \cdot (c(\omega) r^\alpha)^{4/n-2} \]
Thus for $\rho = \frac{n-2}{2 |\alpha| } \cdot c(\omega) \cdot  r^{1 + \alpha \cdot \frac{2}{n-2}}$ we get
\[ \frac{n-2}{2 |\alpha| } \cdot \lambda_H \cdot  c(\omega)^{4 (n-3)/n-2} \cdot a(\omega)^2/\rho^2   = \]
\[ \frac{n-2}{4 (n-1)} scal_{(c(\omega) r^\alpha)^{4/n-2} \cdot g_H}\left(\omega,\left(\frac{2 |\alpha|}{ n-2} \cdot c(\omega)^{-1} \cdot \rho \right)^{1/1 +
\alpha \cdot \frac{2}{n-2}}\right)\] where the right hand side is the new scalar curvature in a point (in the ray $(\omega,t) \in C$) with new distance $\rho$
to $0$. Thus if we rewrite this equation after arc-length reparametrization  (i.e. choosing new coordinates) we get
\[ \frac{4 (n-1)}{2 |\alpha| } \cdot \lambda_H \cdot c(\omega)^{4 (n-3)/n-2} \cdot a(\omega)^2/\rho^2   =
scal_{(c(\omega) r^\alpha)^{4/n-2} \cdot g_H}\left(\omega,\rho \right)\] \qed

When we look at the expression for the scalar curvature we notice that $scal \ge 0$ but it is zero in those places where $|A| = 0$. Using the conformal deformation tools of [L3] sec.2-3 we can
slightly modify the geometry: as in the discussion of funnel sets in sec. 2 above we know that the zero set of $|A|$ is surrounded by a region with lower boundedly increasing $|A|$. This is true
on $\p B_1(0) \cap G_a \cap C$ for each cone $C$ for any $a > 0$. And the definition of $G_a$ shows that these lower bound a re uniform on the space of all minimal cones. Thus what we do
is simply us the amount (=$\frac{4 (n-1)}{2 |\alpha| } \cdot \lambda_H \cdot c(\omega)^{4 (n-3)/n-2} \cdot a(\omega)^2/\rho^2$) of positive scalar curvature close to $|A|^{-1}(0)$ and
conformally deform $C$ on a  uniformly small neighborhood of $|A|^{-1}(0)$ to a new one shifting positive  scalar curvature to
$|A|^{-1}(0)$ decreasing along the border of the neighborhood keeping the scalar curvature positive everywhere and respecting the cone shape. \\

\begin{lemma} \quad   Let  $(C,\tilde g)$ be an abstract cone with  $scal_{(C,\tilde g)} = \frac{4 (n-1)}{2 |\alpha| } \cdot \lambda_H \cdot
c(\omega)^{4 (n-3)/n-2} \cdot a(\omega)^2/\rho^2$. Then we can conformally deform the metric $\tilde g$ to some cone metric $\tilde g^\ast$ with $scal_{\tilde
g^\ast}(\omega,\rho) \ge \iota_H/\rho^2$ on $G_1$ for some $ \iota_H > 0$ which is independent of the singular cone $C  \in \overline {\cal T}_H$.\\ After some additional
reparametrization the metric $\tilde g^\ast$ is again of
the same form  $\tilde c(\omega)^{4/n-2}  \cdot g_{\R} + r^2 \cdot g_{\p B_1(0) \cap C}$\\
\end{lemma}

(Acc. [L3] the scalar curvature of an overlapping collection of such deformations from different tangent cones can still estimated (from below) for upper bounded intersection numbers
(Besicovitch coverings)).

\vspace{1cm}

\setcounter{section}{8}
\renewcommand{\thesubsection}{\thesection}
\subsection{Bridges between $H$ and its tangent cones}

\bigskip

Now we want to prove that $\wp_H$ induces precisely $\wp_C$ on its tangent cones.\\

We already know that there is (up to multiples) only one positive solution on each tangent cone $C$ which is Perron and that arises as a limit of Dirichlet eigenfunctions on $G_a
\subset C$ for
$a \ra \infty$ resp. for an exhausting sequence of regular subcones.  $\wp_H$ can be characterized in the same way  and we want to exploit this fact to prove that $\wp_H$ induces $\wp_C$.\\
Although $w_{k \cdot a}$ induces precisely $w_a$ on each tangent cone for $k \ra \infty$ the fact that $\wp_H$ induces precisely $\wp_C$ is not an immediate consequence since
we can only estimate how fast $w_a$ converges to $\wp_C$ on $C$ whereas a priori this convergence could be that slow on $H$ such that $\wp_H$ induces a function with a higher
order singularity in $0$ than that of $\wp_C$. However a combination with the freezing effect described in and before \ref{approx-lemma} can be used to derive some \emph{telescope} argument.\\

We start with a simple but important \textbf{observation}:  consider a solution  $f = c(\omega) \cdot (\eta_- \cdot r^{\alpha_-} + \eta_+ \cdot r^{\alpha_+})$, for some
$\eta_\pm \ge 0$ on $C$.
 Now, around $0$, we zoom in, that is, we consider $f(\gamma \cdot x)$  for some  $\gamma < 1$ and compare it with $f(x)$  \[c(\omega) \cdot (\eta_- \cdot (\gamma \cdot r)^{\alpha_-} + \eta_+ \cdot (\gamma \cdot r)^{\alpha_+}) =
 c(\omega) \cdot (\eta_- \cdot \gamma^{\alpha_-} \cdot r^{\alpha_-} + \eta_+ \cdot \gamma^{\alpha_+} \cdot r^{\alpha_+})\]
 Compare the ratio of the coefficients for $f(x)$ and $f(\gamma \cdot x)$ (for, say, $\eta_+ > 0$):
 \[\mbox{ since }\; \alpha_- < \alpha_+, \;\;\mbox{ we get } \;  \eta_- /\eta_+ <\gamma^{\alpha_-} \cdot \eta_- /\gamma^{\alpha_+} \cdot \eta_+\]
In other words: For  any linear combination of the two typical solutions,  zooming around $0$ and rewriting the solution in the scaled picture we find the coefficient of the Perron
solution decreases rapidly (to the power of $\alpha_- -  \alpha_+$) relative to the  coefficient of the other generator.\\
Thus normalizing local $L^2$-norms we get for any positive
solution $f = c \cdot \wp + g$, $c > 0$ where $g$ has a lower asymptotic growth rate near infinity than $\wp_C$: $f(\gamma \cdot x)$ converges  for $\gamma \ra 0$ $C^k$-compactly to $\wp_C$.\\
Thus, and we call this a \emph{Perron recovery process} ,if we follow any positive solution on a cone that contains a non vanishing contribution from the Perron solution we can  run along any ray
in $C$ to infinity and observe that (modulo scaling) the function eventually approaches the Perron solution.
Finally note that $w_a$ has the same asymptotic growth rate near infinity as $\wp_C$ which characterizes $\wp_C$ as the Perron solution on $C$ thus this recovery process can clearly be applied to $w_a$. Now can prove\\

\begin{proposition} \quad $\wp_H$ induces precisely $\wp_C$ on any tangent cone $C$.\\
\end{proposition}
\textbf{Proof}\quad We choose a $p \in \Sigma$ and recall the freezing property \ref{approx-lemma} and \ref{aprr}: For any $\delta > 0$ and any triple $R \gg 1 \gg r \gg a >0$ we can find a
small $\eta_{\delta , R , r, a}> 0$ such that for {\it every} $\eta \in (0, \eta_{\delta , R , r, a})$ there is a tangent cone $C^\eta_p$ such that the corresponding part of $\eta^{-2} \cdot H$
can be written as the graph of a function
$g_\eta$ over $C^\eta_p\cap B_R(0)\setminus (B_r (p)\cup V_{a}(\sigma^\eta_p))$ such that $ |g_\eta|_{C^k} < \delta$.\\

Now we fix a tangent cone $C$. By definition this cone appears for a sequence depending on $C$ $r_n \ra 0$, $n \ra \infty$ as a $\delta$-approximating (in the sense detailed above) tangent
cone for $B_{r_n}(p)$
scaled by $(r_n)^{-2}$.  And this also holds for all intermediately appearing tangent cones different form $C$.\\

Fixing the radius $r_1$ we know that for $a \ra \infty$ $w_a$ on $r_1^{-2} \cdot H$ in the part of $H$ identified with $C \cap B_R(0)\setminus (B_r (p)\cup V_{a}(\sigma^\eta_p))$ by
$\delta$-approximation  $C^k$-converges to $\wp_H$. Now we can choose $a_n \ra \infty$ such that this approximation $\wp_H$ by
 $w_{a_n}$ but also by the induced (at this point not yet understood) function on the cone on this fixed part is already very fine and such that $G_a \cap B_{r_n}(p)$ converges for $n \ra \infty$ after rescaling to $G_1 \subset C$.
 Thus in this  region
 we know that  $w_{a_n}$ is approximates $\wp_C$ on $C$ but also  $w_{a_n}$ approximates   $w_{a_n}$ on $H$.\\

 If $C$ were the only, that is a unique, tangent cone for $H$ in $p$ we could now argue directly that on the part of $H$ identified with $C \cap B_R(0)\setminus (B_r (p)\cup V_{a}(\sigma^\eta_p))$
 $\wp_H$ induces $\wp_C$.\\

In general we claim we can find $\delta \ll 1$ and radii $R \gg 1, r \ll 1$ such that we can use the following telescope argument to estimate the growth of the function induces on $C$
by $\wp_H$: the freezing effect also says that the variation of tangent cones slows down while approaching $p \in \Sigma$. That is there are arbitrarily large regions where two
different cones $\delta$-approximate. Thus we start with $w_a$ on $G_a \subset C$ and now run - leaving $p$ - through the family of approximating cones. Starting from $C$
there is such a cone $C_1$ and we use $w_a$ to induce locally a solution on $C_1$ (i.e. transition from $C$ to $H$ and from $H$ to $C_1$). This function differs from a entire
solution i.e. one that is defined all over $C_1$ by an arbitrarily small $C^k$-norm amount (which is uniform for all cones via compactness): namely using the compactness of
$\overline {\cal T}_H \subset\mathcal{C}_{n}$ we can repeat this argument for finer approximations (closer to $p$ and after rescaling) we observe that there is an exhausting
sequence of regular domains in $C_1$ and the sequence of functions defined on these
domains converges to an entire solution on $C_1$.\\
Now we trace this induced function on $C_1$ to infinity, and observe that eventually the asymptotic  growth rate is at most that of $w_a$ on $C$.\\
 The idea is to iterate these cone
transitions until we reach the part fixed above. If we  can manage that the effect of the recovery process  (over)compensates the perturbation of the growth rate during cone transitions we can
clearly argue via growth rates that
$\wp_H$ eventually  approaches $\wp_C$ on the previously fixed part (and to any degree when this presently fixed part is chosen closer to $p$ which is the claim). \\

Thus we have two repelling effects:  the spoiling effect for  the growth rate arising from the transition from $C$ to $C_1$ and the approximation by an entire solution on $C_1$  and
on the other and the recovery process provides an aid to reproduce the desired growth rate.\\

The freezing effect shows that the recovery process eventually (close to $p$) dominates: start with $w_a$ on $C$ wait until far away from $\p G_a$  the growth rate is that of the
Perron solution $\wp_C$ up  to a tiny $\ve > 0$ (say bounded by $c(\omega) \cdot r^{\alpha \pm \ve}$. Consider the induced function on $C_1$ it may have the growth rate of
$\wp_C$ up to $2 \ve$. However it can be assumed to be close to an entire solution induced from $w_a$ on $C$ whose growth rate near infinity in $C_1$ is again that of $\wp_C$.
Thus we wait until for the induced solution on $C_1$ the growth rate is again that of the Perron solution $\wp_C$ up  to $\ve > 0$. \\
The distance between the region where a
tangent cone receives a part of a solution from another cone from the transition and the region where the recovery process reproduced the $\ve$-variation from the growth of
the Perron solution is uniformly upper bounded using the compactness of $\overline {\cal T}_H \subset\mathcal{C}_{n}$. Thus for a sufficiently small starting radius $r_1$ this
process shows that $\wp_H$ eventually approaches $\wp_C$ to any desired degree of accuracy. \qed

As a certainly expects but not yet discussed consequence of this argument we notice that the Perron solutions on all tangent cones in a given point $p \in \Sigma$ have the
same growth rates.

\end{document}